\documentclass[11pt]{article}
\usepackage{amsmath}
\usepackage{amssymb}
\usepackage{amsthm}
\usepackage{url}
\begin{document}
\newtheorem{theorem}{Theorem}
\newtheorem{lemma}[theorem]{Lemma}
\newtheorem{corollary}[theorem]{Corollary}
\newtheorem{definition}[theorem]{Definition}
\newtheorem{example}[theorem]{Example}
\pagenumbering{roman}
\renewcommand{\thetheorem}{\thesection.\arabic{theorem}}
\renewcommand{\thelemma}{\thesection.\arabic{lemma}}
\renewcommand{\thedefinition}{\thesection.\arabic{definition}}
\renewcommand{\theexample}{\thesection.\arabic{example}}
\renewcommand{\theequation}{\thesection.\arabic{equation}}
\newcommand{\mysection}[1]{\section{#1}\setcounter{equation}{0}
\setcounter{theorem}{0} \setcounter{lemma}{0}
\setcounter{definition}{0}}
\newcommand{\mrm}{\mathrm}
\newcommand{\beq}{\begin{equation}}
\newcommand{\eeq}{\end{equation}}
\def\dd{\boldsymbol{d}}
\newcommand{\ben}{\begin{enumerate}}
\newcommand{\een}{\end{enumerate}}

\newcommand{\beqa}{\begin{eqnarray}}
\newcommand{\eeqa}{\end{eqnarray}}

\newcommand {\non}{\nonumber}
\newcommand{\C}{\mbox{$\mathbb{C}$}}
\title
{\bf Acceleration of Convergence  of Some Infinite Sequences $\boldsymbol{\{A_n\}}$ Whose  Asymptotic Expansions  Involve Fractional Powers of $\boldsymbol{n}$ via the ${\tilde{\dd}}^{(m)}$~Transformation}
\author
{Avram Sidi\\
Computer Science Department\\
Technion - Israel Institute of Technology\\ Haifa 32000, Israel\\~ \\
E-mail:\ \ \url{asidi@cs.technion.ac.il}\\
URL:\ \ \url{http://www.cs.technion.ac.il/~asidi}}
\date{Appeared in: \ {\em Numerical Algorithms}, 85:1409--1445, 2020}
\maketitle \thispagestyle{empty}
\newpage
\begin{abstract}
In this paper, we discuss the application of the author's $\tilde{d}^{(m)}$ transformation to accelerate the convergence of infinite series $\sum^\infty_{n=1}a_n$ when the terms $a_n$ have asymptotic expansions that can be expressed in  the form
$$ a_n\sim(n!)^{s/m}\exp\left[\sum^{m}_{i=0}q_in^{i/m}\right]\sum^\infty_{i=0}w_i n^{\gamma-i/m}\quad\text{as $n\to\infty$},\quad s\ \text{integer.}$$
We discuss the implementation  of the $\tilde{d}^{(m)}$ transformation via the recursive W-algorithm of the author. We show how to apply this transformation and how to assess in a reliable way the accuracies  of the  approximations it produces, whether the series converge or they diverge.
We classify  the different cases that exhibit unique numerical stability issues  in floating-point arithmetic.
We show that the $\tilde{d}^{(m)}$ transformation can also be used efficiently to accelerate the convergence of  infinite products
 $\prod^\infty_{n=1}(1+v_n)$, where
 $v_n\sim \sum^\infty_{i=0}e_in^{-t/m-i/m}$ as $n\to\infty$,\  $t\geq m+1$ an integer.
Finally, we give several numerical examples that attest  the high efficiency of the
$\tilde{d}^{(m)}$ transformation for the different cases.

\end{abstract}
\vspace{3cm} \noindent {\bf Mathematics Subject Classification
2010:}  40A05, 40A25, 40A20, 65B05, 65Bl0.

\vspace{1cm} \noindent {\bf Keywords and expressions:} Acceleration of convergence,
extrapolation, infinite series, infinite products,
asymptotic expansions, fractional powers, $\tilde{d}^{(m)}$ transformation, W-algorithm.

 \thispagestyle{empty}
\newpage
\pagenumbering{arabic}
\section{Introduction}\label{se1}
\setcounter{equation}{0}
\setcounter{theorem}{0}

The summation of infinite series $\sum^\infty_{n=1}a_n$, where the terms $a_n$ are in general complex and have  asymptotic expansions (as $n\to\infty$) involving  powers of $n^{-1/m}$ for  positive integers  $m$, has been
 of some interest. Due to their complex analytical nature, however, the rigorous study of such series has been the subject of very few works. See Birkhoff \cite{Birkhoff:1930:FTI}  and
  Birkhoff and Trjitzinsky \cite{Birkhoff:1932:ATS}. For a brief summary of these works, see  Wimp \cite{Wimp:1974:SSW}, \cite[Section 1.7]{Wimp:1981:STA}.

  In this work, we deal with those infinite series $\sum^\infty_{n=1}a_n$, whether convergent or divergent, for which  $\{a_n\}$  belong to a class of sequences denoted $\tilde{\bf b}^{(m)}$. These series were first   studied in detail   in Sidi \cite[Section 6.6]{Sidi:2003:PEM}, where an extrapolation method denoted the $\tilde{d}^{(m)}$ transformation  to accelerate their convergence (actually, to accelerate the convergence of the sequence $\{A_n\}$ of the partial sums $A_n=\sum^n_{k=1}a_k$, $n=1,2,\ldots$)  was also developed. This transformation is
  very effective also when these series diverge; in such cases, it produces approximations to the antilimits of the series treated.
Practically speaking, a sequence $\{a_n\}$ is in $\tilde{\bf b}^{(m)}$, $m\geq1$ being an integer, if
$a_n$ has an asymptotic expansion that {\em can be expressed} in the form
\beq \label{eq:ff7r}
a_n\sim[\Gamma(n)]^{s/m}\exp\left[Q(n)\right]\sum^\infty_{i=0}w_i n^{\gamma-i/m}\quad\text{as $n\to\infty$},\eeq
where
\begin{enumerate}
\item$\Gamma(z)$ is the gamma function,
\item
$s$ is an arbitrary integer, positive, negative, or zero,
\item
$Q(n)$ is either identically zero or  is a polynomial of degree at most $m$ in $n^{1/m}$, expressed as in
\beq\label{eq:ff8r}
Q(n)=\sum^{m-1}_{i=0}\theta_in^{1-i/m},\eeq
$\theta_0,\theta_1,\ldots,\theta_{m-1}$ being  real or complex constants,\footnote{Clearly, $Q(n)\equiv0$ takes place when $\theta_i=0$, $i=0,1,\ldots,m-1.$}
\item
$\gamma$ is an arbitrary real or complex number.
\end{enumerate}

In the special case of $m=1$,  either $Q(n)=\theta_0 n$ with $\theta_0\neq0$  or $Q(n)\equiv0$, and  \eqref{eq:ff7r} assumes the form
\beq \label{eq:ff7w}
a_n\sim[\Gamma(n)]^{s}\zeta^n\sum^\infty_{i=0}w_i n^{\gamma-i}\quad\text{as $n\to\infty$},\eeq with (i)\,$\zeta=1$ if $Q(n)\equiv0$ and (ii)\,$\zeta=e^{\theta_0}\neq1$ if $Q(n)=\theta_0 n$ with $\theta_0\neq0$.
Here we note that
the class $\tilde{\bf b}^{(1)}$ is simply the class denoted ${\bf b}^{(1)}$,  which is a special case and the simplest prototype  of
 the collection of sequence classes ${\bf b}^{(m)}$,  $m=1,2,\ldots,$ originally introduced in Levin  and Sidi \cite{Levin:1981:TNC} and   studied extensively in   Sidi \cite[Chapter 6]{Sidi:2003:PEM}.\footnote{The convergence of infinite series $\sum^\infty_{n=1}a_n$ with $\{a_n\}\in {\bf b}^{(m)}$,   $m$ being arbitrary, can be accelerated
 efficiently by using the $d^{(m)}$ transformation of Levin  and Sidi \cite{Levin:1981:TNC}, which can be implemented very economically via the recursive W$^{(m)}$-algorithm of Ford and Sidi \cite{Ford:1987:AGR}. All this  is  studied in detail also  in   Sidi \cite[Chapters  6 and 7]{Sidi:2003:PEM}.}
In this connection, we mention that the $t$, $u$, and $v$ transformations of Levin
\cite{Levin:1973:DNT}  and the $d^{(1)}$ transformation of Levin
 and Sidi \cite{Levin:1981:TNC} are very effective convergence acceleration methods
 for infinite series $\sum^\infty_{n=1}a_n$ with $\{a_n\}\in{\bf b}^{(1)}$.

In this work, we shall deal with the class  $\tilde{\bf b}^{(m)}$, $m\geq1$ being arbitrary. We shall use the notation of \cite[Section 6.6]{Sidi:2003:PEM} throughout. Comparing \eqref{eq:ff7r}--\eqref{eq:ff8r} with \eqref{eq:ff7w}, and judging also from
Theorem \ref{th:ff3}, we realize that sequences in  $\tilde{\bf b}^{(m)}$ with $m\geq2$  have a  richer and
more interesting mathematical structure than those in $\tilde{\bf b}^{(1)}={\bf b}^{(1)}$. As  will also be  clear from the numerical examples in Section \ref{se5}, depending on whether $a_n$ in \eqref{eq:ff7r} is such that
\begin{itemize}
\item [(i)]
$s=0$ and $Q(n)\equiv0$ and $\gamma\neq -1+i/m$, $i=0,1,\ldots,$ or
\item [(ii)]
$s=0$ and $Q(n)\not\equiv0$, with $\theta_0\neq0$ and $\gamma$ is arbitrary, or
\item [(iii)]
$s=0$ and $Q(n)\not\equiv0$, with $\theta_0=\cdots=\theta_{r-1}=0$ and $\theta_r\neq0$ for some $r\in\{1,\ldots,m-1\}$, and $\gamma$ is arbitrary,  or
\item [(iv)]
$s\neq0$ ($s<0$ or $s>0$) and $Q(n)$ is arbitrary [$Q(n)\equiv0$ or
 $Q(n)\not\equiv0$], and $\gamma$ is arbitrary, or
 \item [(v)]
 $a_n$ is as in any one of the cases (i)--(iv) (with real $\theta_0$), multiplied by $(-1)^n$,
\end{itemize}
 the series $\sum^\infty_{n=1}a_n$
 exhibit  different convergence and numerical stability properties when convergence acceleration methods are applied to them in finite-precision (floating-point) arithmetic.  In addition, the series $\sum^\infty_{n=1}a_n$ may converge or diverge.

The contents of this paper are arranged as follows:
In the next section, we  summarize the asymptotic properties of sequences $\{a_n\}$ in  $\tilde{\bf b}^{(m)}$ for arbitrary $m$.  In Section~\ref{se3},  (i)\,we recall   the  $\tilde{d}^{(m)}$ transformation, (ii)\,we recall the issue of assessing the numerical stability of the approximations generated by it, (iii)\,we  recall the  W-algorithm of Sidi \cite{Sidi:1982:ASC} as it is used for implementing the $\tilde{d}^{(m)}$ transformation, and  (iv)\,we discuss how the  W-algorithm can be  extended for assessing in a very simple way  the numerical stability of the approximations generated by the $\tilde{d}^{(m)}$ transformation simultaneously with their computation in finite-precision arithmetic.
In Section \ref{se4},  we  illustrate Theorem \ref{th:ff4}, which concerns the asymptotic behavior of the partial sums $A_n=\sum^n_{k=1}a_k$ as $n\to\infty$, on the basis of which the
$\tilde{d}^{(m)}$ transformation is developed, with some instructive  examples. In Section \ref{se5},  we  illustrate with numerical examples of varying nature the remarkable effectiveness of the  $\tilde{d}^{(m)}$ transformation on the series $\sum^\infty_{n=1}a_n$, where  $\{a_n\}\in\tilde{\bf b}^{(m)}$, whether these converge or diverge. We also show how the $\tilde{d}^{(m)}$ transformation can be tuned for best numerical results. In Section \ref{se6}, we consider the use of the
$\tilde{d}^{(m)}$
transformation for computing   some infinite products $\prod^\infty_{n=1}(1+v_n)$, where $\{v_n\}\in \tilde{\bf A}_0^{(-t/m,m)}$, that is,
\beq \label{eqzx1}v_n\sim\sum^\infty_{i=0}w_in^{-t/m-i/m},\quad t\geq m+1 \ \ \text{an integer.}\eeq
 We study the asymptotic behavior of the partial products $A_n=\prod^n_{k=1}(1+v_k)$ as $n\to\infty$ and conclude  that   the $\tilde{d}^{(m)}$
transformation can be applied very efficiently to accelerate the convergence of the sequence of
the partial products. In Section \ref{se7}, we give numerical examples that illustrate the efficiency of the $\tilde{d}^{(m)}$ transformation on such infinite products.

Presently, there is no  numerical experience with
the issue of convergence acceleration of the infinite series described above in their most general form, that is, with arbitrary $m$, $s$, $\gamma$, and $Q(n)$.   So far, the acceleration of the convergence of only  a  subset of such series, for which   $s=0$ and $Q(n)\equiv0$ and $\sum^\infty_{n=1}a_n$ is convergent,    has been considered in the literature; thus
\beq \label{eqzx2}a_n\sim\sum^\infty_{i=0}w_i n^{\gamma-i/m}\quad\text{as $n\to\infty$}\quad\text{and}\quad  \Re\gamma <-1,\eeq in this subset:  Sablonni\`{e}re
 \cite{Sablonniere:1992:ABI} has studied  the application of (i)\,the  iterated modified     $\Delta^2$-process  and
 (ii)\,the iterated $\theta_2$-algorithm of Brezinski \cite{Brezinski:1975:GTS},  to the cases in which $m=1,2$ only.
   Van Tuyl \cite{VanTuyl:1976:AMA}, \cite{VanTuyl:1994:ACF} has studied the application of (i)\,the  iterated  modified     $\Delta^2$-process, (ii)\,the  iterated
   transformation of Lubkin \cite{Lubkin:1952:MSI}, (iii)\,the $\theta$-algorithm of Brezinski \cite{Brezinski:1975:GTS}, (iv)\,a generalization of the $\rho$-algorithm of Wynn \cite{Wynn:1956:PTN}, (v)\,the $u$ and $v$ transformations of Levin \cite{Levin:1973:DNT}, (vi)\,a generalization of the
   Neville table,  and (vii)\,the $d^{(m)}$ transformation of Levin and Sidi \cite{Levin:1981:TNC}.
   The numerical results of \cite{VanTuyl:1976:AMA} show that, with the exception of the $u$ and $v$ transformations, which are effective only when $m=1$, the rest of the transformations are effective accelerators  for all $m\geq1$.
 (Note that the iterated $\theta_2$-algorithm and iterated
     Lubkin transformation are identical.)

   The modified  $\Delta^2$-process is due to Drummond \cite{Drummond:1976:SCT} (see also
    Brezinski and Redivo-Zaglia \cite{Brezinski:2010:EDP}),  while the generalized $\rho$-algorithm and
   the generalized Neville table are given in  Van Tuyl \cite{VanTuyl:1994:ACF}.
   For the  $\Delta^2$-process, which is due to Aitken \cite{Aitken:1926:BNS}, see Stoer and Bulirsch
   \cite[Chapter 5]{Stoer:2002:INA} and    Sidi \cite[Chapter 15]{Sidi:2003:PEM},
 for example.
    For  discussions of  the methods mentioned above, see  also \cite[Chapters 6, 15, 19, 20]{Sidi:2003:PEM}.

    We note that to apply the  modified     $\Delta^2$-process,
     the generalized $\rho$-algorithm, and the generalized Neville table, we need to know $\gamma$ in \eqref{eqzx2}. This is not the case when applying the iterated transformation of Lubkin,  the $\theta$-algorithm,   the  $d^{(m)}$ transformation, and the $\tilde{d}^{(m)}$ transformation.

     Before proceeding further, we would like to emphasize that the
  $\tilde{d}^{(m)}$ transformation can be formulated such that   it will be applicable without any modification and  with success  to {\em all} infinite series $\sum^\infty_{n=1}a_n$ where  $\{a_n\}\in\tilde{\bf b}^{(m)}$,  with {\em arbitrary} $s$, $Q(n)$, and $\gamma$, which do {\em not} have to be known. This is a very important feature of the $\tilde{d}^{(m)}$ transformation and of this work.

Finally, we mention that the works \cite{Sablonniere:1992:ABI} and \cite{VanTuyl:1994:ACF} deal only with the convergence issue of the transformations discussed in them; they do not consider the important issue of numerical stability when using  floating-point (finite-precision) arithmetic.\footnote{Note that most of the methods mentioned above suffer from lack of numerical stability when applied to infinite series
$\sum^\infty_{n=1}a_n$ with $a_n$ behaving as in \eqref{eqzx2}. In addition, there is no  reliable way to assess the  floating-point accuracies of the approximations they produce.}
In our treatment of the $\tilde{d}^{(m)}$ transformation in Section \ref{se3} of this work, we emphasize this issue as follows: (i)\,we devise  reliable zero-cost procedures for monitoring  the numerical stability and predicting the maximum accuracy of the approximations produced at the time these are being computed and (ii)\,we overcome numerical instabilities  by applying the  $\tilde{d}^{(m)}$ transformation to  properly sampled subsequences of  the sequences $\{A_n\}$ of partial sums $A_n=\sum^n_{k=1}a_k$ via  {\em arithmetic progression sampling (APS)} or  {\em geometric progression sampling (GPS)}--- two automatic sampling procedures originally proposed in Ford and Sidi \cite{Ford:1987:AGR}---that have been shown to be very effective.   These are two additional important features of this work that differentiate it from all previous works.

\section{Preliminaries} \label{se2}
\setcounter{equation}{0}
\setcounter{theorem}{0}
\subsection{The function class $\tilde{\mathbf{A}}^{(\gamma,m)}_0$}
 We begin with the following definition:

  \begin{definition} [{\bf\cite{Sidi:2003:PEM}, Definition 6.6.1}]\label{def:ff1}
A function $\alpha(x)$ defined for all large $x$
is in the set $\tilde{\mathbf{A}}^{(\gamma,m)}_0$,
$m$ being a positive integer,
if it has a Poincar\'{e}-type asymptotic expansion of
the form
\begin{equation} \label{eq:ff1}
\alpha(x) \sim  \sum^\infty_{i=0}\alpha_i
x^{\gamma -i/m}\ \ \mbox{as}\ x\to \infty.
\end{equation}
In addition, if $\alpha_0\neq 0$ in \eqref{eq:ff1}, then
$\alpha(x)$ is said to belong to
$\tilde{\mathbf{A}}^{(\gamma,m)}_0$ strictly. Here $\gamma$ is
complex in general.\footnote{Clearly, if $\alpha\in\tilde{\mathbf{A}}^{(\gamma,m)}_0$, then
$\alpha(x)=x^\gamma\beta(x)$, where $\beta\in\tilde{\mathbf{A}}^{(0,m)}_0$.}
\end{definition}

Before going on, we state some properties of the sets $\tilde{\mathbf{A}}^{(\gamma,m)}_0$, whose verification we leave to the reader. We make repeated use of these properties  in Sections \ref{se4} and \ref{se6}.
\begin{enumerate}
\item \label{re:f1-first}
$\tilde{\mathbf{A}}^{(\gamma,m)}_0 \supset\tilde{\mathbf{A}}^{(\gamma-1/m,m)}_0 \supset
\tilde{\mathbf{A}}^{(\gamma-2/m,m)}_0 \supset\cdots,$ so
that if $\alpha \in \tilde{\mathbf{A}}^{(\gamma,m)}_0$, then, for any positive integer $k$,
$\alpha \in \tilde{\mathbf{A}}^{(\gamma+k/m,m)}_0$ but not strictly.
Conversely, if
$\alpha \in \tilde{\mathbf{A}}^{(\delta,m)}_0$ but not strictly, then
$\alpha \in \tilde{\mathbf{A}}^{(\delta-k/m,m)}_0$ strictly for a unique positive integer $k$.
\item If $\alpha \in \tilde{\mathbf{A}}^{(\gamma,m)}_0$ strictly, then $\alpha \notin
\tilde{\mathbf{A}}^{(\gamma -1/m,m)}_0$.
\item If $\alpha \in \tilde{\mathbf{A}}^{(\gamma,m)}_0$ strictly, and $\beta(x) = \alpha
(cx+d)$ for some arbitrary constants $c>0$ and $d$, then $\beta \in
\tilde{\mathbf{A}}^{(\gamma,m)}_0$ strictly as well.
\item If $\alpha,\beta \in \tilde{\mathbf{A}}^{(\gamma,m)}_0$, then $\alpha \pm \beta \in
\tilde{\mathbf{A}}^{(\gamma,m)}_0$ as well.
(This implies that the
zero function is included in $\tilde{\mathbf{A}}^{(\gamma,m)}_0$.)
If $\alpha \in\tilde{\mathbf{A}}^{(\gamma,m)}_0$
and $\beta \in \tilde{\mathbf{A}} ^{(\gamma+k/m,m)}_0$ strictly for
some positive integer $k$, then $\alpha \pm \beta \in\tilde{\mathbf{A}}^{(\gamma+k/m,m)}_0$
strictly.

\item If $\alpha \in\tilde{\mathbf{A}}^{(\gamma,m)}_0$ and $\beta \in\tilde{\mathbf{A}}^{(\delta,m)}_0$,
then $\alpha \beta \in\tilde{\mathbf{A}}^{(\gamma+\delta,m)}_0$;  if, in addition,
$\beta \in\tilde{\mathbf{A}}^{(\delta,m)}_0$ strictly, then $\alpha / \beta \in
\tilde{\mathbf{A}}^{(\gamma - \delta,m)}_0$.
\item If $\alpha \in\tilde{\mathbf{A}}^{(\gamma,m)}_0$ strictly, such that $\alpha(x)>0$
for all large $x$, and we define $\theta(x)=
{[\alpha (x)]}^\xi$, then $\theta \in\tilde{\mathbf{A}}^{(\gamma \xi,m)}_0$
strictly.
\item If $\alpha \in\tilde{\mathbf{A}}^{(\gamma,m)}_0$ strictly and $\beta \in
\tilde{\mathbf{A}}^{(k,1)}_0$ strictly for some positive integer $k$, such that
$\beta(x)>0$ for all large $x>0$,
and we define $\theta(x)=\alpha(\beta(x))$,
then $\theta\in \tilde{\mathbf{A}}^{(k\gamma,m)}_0$ strictly.

\item
If $\alpha \in\tilde{\mathbf{A}}^{(\gamma,m)}_0$ (strictly), and $\beta(x)=\alpha(x+d)-
\alpha(x)$ for an arbitrary constant $d\neq0$, then $\beta\in\tilde{\mathbf{A}}^{(\gamma-1,m)}_0$
(strictly) when $\gamma\neq0$. If $\alpha \in\tilde{\mathbf{A}}^{(0,m)}_0$, then
$\beta\in\tilde{\mathbf{A}}^{(-1-1/m,m)}_0$.
\end{enumerate}

Note that if $a_n=\alpha(n)$, $n=1,2,\ldots,$ where $\alpha\in\tilde{\mathbf{A}}^{(\gamma,m)}_0$, then $a_n$ is as in \eqref{eqzx2}.
Such sequences $\{a_n\}$  are therefore in the class $\tilde{\mathbf{b}}^{(m)}$.

The following theorem summarizes the summation properties of functions  in
$\tilde{\mathbf{A}}^{(\gamma,m)}_0$. It is also useful in proving Theorem \ref{th:ff2}.
\begin{theorem} [{\bf\cite{Sidi:2003:PEM}, Theorem 6.6.2}] \label{th:ff1}
Let $g\in \tilde{\mathbf{A}}
^{(\gamma,m)}_0$ strictly for some $\gamma$ with
$g(x)\sim \sum^\infty_{i=0}
g_ix^{\gamma-i/m}$ as $x\to \infty$, and  define
$G(n)=\sum^{n-1}_{r=1}g(r)$. Then
\begin{equation}\label{eq:ff2}
G(n)=b+c \log n+\tilde{G}(n),
\end{equation}
where $b$ and $c$ are constants and $\tilde{G}\in
\tilde{\mathbf{A}}^{(\gamma+1,m)}_0$.

\begin{enumerate}
\item
If $\gamma \neq -1$, then $ \tilde{G}\in \tilde{\mathbf{A}}_0
^{(\gamma+1,m)}$ strictly, while
$\tilde{G}\in \tilde{\mathbf{A}}^{(-1/m,m)}_0$ if $\gamma=-1$.

\item
If  $\gamma+1\neq i/m,\ i=0,1,\ldots,$ then  $c=0$, and either
(i)\,$b$ is the limit  of $G(n)$ as $n\to
\infty$ if $\Re\gamma+1<0$, or (ii)\,$b$ is the  antilimit of $G(n)$ as $n\to
\infty$ if $\Re\gamma+1\geq0.$
\item
If $\gamma+1=k/m$ for some integer
$k\geq 0,$ then $c=g_k$.
\end{enumerate}
Finally,
\begin{equation}\label{eq:ff3}
\tilde{G}(n)=\sum^{m-1}_{\substack{i=0 \\{\gamma-i/m\neq -1}}}
\frac{g_i}{\gamma-i/m+1}n^{\gamma-i/m+1}+O(n^\gamma)\ \
\mbox{as}\ n\to \infty.
\end{equation}
\end{theorem}

Before ending this section, we also note that  the sets $\tilde{\mathbf{A}}^{(\gamma,m)}_0$ are most important building blocks of sequences $\{a_n\}$  in the  class $\tilde{\bf b}^{(m)}$,  to which we turn next.

\subsection{The sequence class $\tilde{\mathbf{b}}^{(m)}$}

With the classes $\tilde{\mathbf{A}}^{(\gamma,m)}_0$ already
defined, we now go on to define the sequence class
$\tilde{\mathbf{b}}^{(m)}$.

\begin{definition} [{\bf\cite{Sidi:2003:PEM}, Definition 6.6.3}]\label{def:ff2}

A sequence $\{a_n\}$ belongs to the class
$\tilde{\mathbf{b}}^{(m)}$ if it satisfies a linear
homogeneous difference equation of first order of the form
$a_n=p(n)\Delta a_n$ with $p\in \tilde{\mathbf{A}}^
{(q/m,m)}_0$ for some integer $q \leq m$.
Here $\Delta a_n=a_{n+1}-a_n$, $n=1,2,\ldots.$
\end{definition}

We begin with the following general result:

\begin{theorem} [{\bf \cite{Sidi:2003:PEM}, Theorem 6.6.4}]\label{th:ff2}
\begin{enumerate}
\item [(i)] Let $a_{n+1}=c(n)a_n$ such that $c\in \tilde{\mathbf{A}}
^{(\mu,m)}_0$ strictly with $\mu$ in general complex.
Then $a_n$ is of the form
\begin{equation}\label{eq:ff7}
a_n=[\Gamma(n)]^\mu \exp \left[Q(n)\right]n^\gamma w(n),
\end{equation}
where $\Gamma(n)$ is the gamma function and $Q(n)$ is a polynomial of degree at most $m$ in $n^{-1/m}$ which we choose to write in the form
\begin{equation}\label{eq:ff8}
Q(n)=\sum^{m-1}_{i=0}\theta_in^{1-i/m},\eeq
and
\beq \label{eq:ff8a}
w\in \tilde{\mathbf{A}}^{(0,m)}_0\ \ \mbox{strictly.}
\end{equation}
Given that $c(n)\sim \sum^\infty_{i=0}c_in^{\mu-i/m}\ \
\mbox{as}\ \ n\to \infty,$ with $c_0\neq0$,
we have
\begin{equation}\label{eq:ff9}
e^{\theta_0}=c_0;\ \ \theta_i=\frac{\epsilon_i}{1-i/m},\ \
i=1,\ldots,m-1;\ \gamma=\epsilon_m,
\end{equation}
where the $\epsilon_i$ are determined by $c_0,c_1,\ldots,c_m$  via
\begin{equation}\label{eq:ff10}
\sum^m_{s=1}\frac{{(-1)}^{s+1}}{s}{\biggl( \sum^m_{i=1}
\frac{c_i}{c_0}z^i \biggr )}^s=\sum^m_{i=1}\epsilon_iz^i
+O(z^{m+1})\ \ \mbox{as}\ z\to 0.
\end{equation}
(Note that $\theta_0=0$ when $c_0=1$.)
\item [(ii)] The converse is also true, that is, if $a_n$ is as in
\eqref{eq:ff7}--\eqref{eq:ff8a}, then $a_{n+1}=c(n)a_n$
with $c\in \tilde{\mathbf{A}}^{(\mu,m)}_0$ strictly.\\
\item [(iii)] Finally, (a)\ \ $\theta_1=\cdots=\theta_{m-1}=0$ if and
only if $c_1=\cdots=c_{m-1}=0,$ and (b)\ \  $\theta_1=\cdots=
\theta_{r-1}=0$ and $\theta_r\neq0$ if and only if
$c_1=\cdots=c_{r-1}=0$ and $c_r\neq0,\ r\in
\{1,\ldots,m-1\}$.
\end{enumerate}
\end{theorem}
\medskip
\noindent{\bf Remark.}  Note that we can express \eqref{eq:ff7} also in the form
$$ a_n=[\Gamma(n)]^\mu \exp \left[\hat{Q}(n)\right]n^\gamma w(n)\, \zeta^n,$$
where
$$
\hat{Q}(n)=\sum^{m-1}_{i=1}\theta_in^{1-i/m},\quad \zeta=c_0=e^{\theta_0}.$$
Of course, $\zeta=1$ when $c_0=1$ and $\zeta\neq1$ when $c_0\neq1$.
\bigskip

The next theorem gives necessary and sufficient conditions for a
sequence $\{a_n\}$ to be in $\tilde{\mathbf{b}}^{(m)}$. In this sense,  it
is a characterization theorem for sequences in
$\tilde{\mathbf{b}}^{(m)}$. Theorem \ref{th:ff2} becomes useful in the
proof.

\begin{theorem}[{\bf \cite{Sidi:2003:PEM}, Theorem 6.6.5}] \label{th:ff3} A sequence $\{a_n\}$ is in
$\tilde{\mathbf{b}}^{(m)}$ if and only if its members satisfy
$a_{n+1}=c(n)a_n$ with $c\in \tilde{\mathbf{A}}^{(s/m,m)}_0$ for some
 integer $s$ and $c(n)\neq 1+O(n^{-1-1/m})$ as
$n\to \infty$.\footnote{Thus, $\{a_n\}\not\in\tilde{\mathbf{b}}^{(m)}$ if
$c(n)=a_{n+1}/a_n$ is  in $\tilde{\mathbf{A}}^{(0,m)}_0$
 and has an asymptotic expansion of the form
 $c(n)\sim 1+\sum^\infty_{i=m+1}c_i n^{-i/m}$ as $n\to\infty$.} Therefore,  $a_n$ is as in \eqref{eq:ff7}--\eqref{eq:ff8a} with $\mu=s/m$, and this implies that
$a_n=p(n)\Delta a_n$ with $p\in\tilde{\mathbf{A}}^{(\sigma,m)}_0$
strictly, where $\sigma=q/m$ and $q$ is an integer $\leq m$.
With $c(n)\sim\sum^\infty_{i=0}c_in^{s/m-i/m}$ as $n\to\infty$, $c_0\neq0$,
 we have the following specific cases:
\begin{enumerate}
\item \label{re1}
When $s=0$,  $c_0=1$, $c_1=\cdots=c_{m-1}=0$, and $c_m\neq0$, which holds necessarily,  we have $\sigma=1$ or $q=m$.\\
 In this case,
$a_n=n^\gamma w(n)$ with $\gamma=c_m\neq0$. Hence $a_n=h(n)$, $h\in \tilde{\mathbf{A}}^{(\gamma,m)}_0$, $\gamma\neq0$.
\item \label{re2}
When $s=0$,  $c_0=1$, $c_1=\cdots=c_{r-1}=0$, and $c_r\neq0$, $r<m$, we have $\sigma=r/m$ or $q=r$.\\
In this case,
$a_n=\exp[Q(n)] n^\gamma w(n)$,
$Q(n)=\sum^{m-1}_{i=r}\theta_i n^{1-i/m}$.

\item \label{re3}
When $s=0$, $c_0\neq1$, we have $\sigma=0$ or $q=0$.\\
In this case,
$a_n=\exp[Q(n)] n^\gamma w(n)$, $Q(n)=\sum^{m-1}_{i=0}\theta_i n^{1-i/m}$, $\theta_0\neq0$.
\item \label{re4}
When $s<0$,  we have $\sigma=0$ or $q=0$.\\
In this case,
$a_n=[\Gamma(n)]^{s/m}\exp[Q(n)] n^\gamma w(n)$, $Q(n)=\sum^{m-1}_{i=0}\theta_i n^{1-i/m}$.

\item \label{re5}
When $s>0$,  we have $\sigma=-s/m$ or $q=-s$.\\
In this case,
$a_n=[\Gamma(n)]^{s/m}\exp[Q(n)] n^\gamma w(n)$, $Q(n)=\sum^{m-1}_{i=0}\theta_i n^{1-i/m}$.

 \end{enumerate}
 Of course, $w(n)$ is in $\tilde{\mathbf{A}}^{(0,m)}_0$ in all  cases.
\end{theorem}

 We now  restrict our attention to the cases described in parts \ref{re1}--\ref{re4} of Theorem \ref{th:ff3}, for which the series $\sum^\infty_{n=1}a_n$ (i)\,either converges (ii)\,or diverges but has  an  Abel sum or an  Hadamard finite part that serves as the antilimit of $A_n=\sum^{n}_{k=1}a_k$ as $n\to\infty$.  (In part \ref{re5},   the series $\sum^\infty_{n=1}a_n$  always diverges and has no  Abel sum or Hadamard finite part. It may have a Borel sum, however.)

\begin{enumerate}
\item
In part  \ref{re1}, we assume that $\gamma\neq -1+i/m$, $i=0,1,\ldots,$ as in part of Theorem \ref{th:ff1}. We have two cases to consider:
\begin{itemize}
\item
If $\Re\gamma<-1$,   $\sum^\infty_{n=1}a_n$  converges.
\item
If $\Re\gamma\geq-1$,  $\sum^\infty_{n=1}a_n$ diverges but has an Hadamard finite part that serves as the antilimit of $A_n=\sum^{n}_{k=1}a_k$ as $n\to\infty$.
\end{itemize}
\item
In part \ref{re2}, we assume the following two situations:
\begin{itemize}
\item
$\lim_{n\to\infty}\Re Q(n)=-\infty$
or, equivalently, $\Re\theta_r<0$. In this case, $\sum^\infty_{n=1}a_n$  converges for all $\gamma$. [If $\Re\theta_r>0$, then $\lim_{n\to\infty}\Re Q(n)= +\infty$; therefore, $\sum^\infty_{n=1}a_n$  diverges for all $\gamma$.]
\item
$\Re Q(n)=0$ or, equivalently, $\Re\theta_i=0$, $i=r,\ldots,m-1$.
\begin{itemize}
\item
$\sum^\infty_{n=1}a_n$  converges if $\Re\gamma<-r/m$.
\item
$\sum^\infty_{n=1}a_n$  diverges if $\Re\gamma\geq-r/m$ but has an Abel sum that serves as the antilimit of $A_n=\sum^{n}_{k=1}a_k$ as $n\to\infty$.
\end{itemize}
\end{itemize}
\item
In part \ref{re3}, as in item \ref{re2}, we assume the following two situations:
\begin{itemize}
\item
$\lim_{n\to\infty}\Re Q(n)=-\infty$
or, equivalently, $\Re\theta_0<0$, which is equivalent to $|c_0|<1$. In this case, $\sum^\infty_{n=1}a_n$  converges for all $\gamma$. (If $\Re\theta_0>0$, which is equivalent to $|c_0|>1$, $\sum^\infty_{n=1}a_n$  diverges for all $\gamma$.)
\item
$\Re Q(n)=0$ or, equivalently, $\Re\theta_i=0$, $i=0,1,\ldots,m-1$.
(Note that we now  have $|c_0|=1$, in addition to $c_0\neq1$.) We now have the following cases:
\begin{itemize}
\item
$\sum^\infty_{n=1}a_n$  converges if $\Re\gamma<0$.
\item
$\sum^\infty_{n=1}a_n$  diverges if $\Re\gamma\geq0$ but has an Abel sum that serves as the antilimit of $A_n=\sum^{n}_{k=1}a_k$ as $n\to\infty$.
\end{itemize}
\end{itemize}

\item
In part \ref{re4}, we do not assume anything in addition to what is there. In this case, $\sum^\infty_{n=1}a_n$  converges for all $\gamma$.
 \end{enumerate}

\noindent{\bf Remark:} Note that in all the cases considered above, we have $\sigma=q/m$, with $q\in\{0,1,\ldots,m\}$.
\medskip

Theorem  \ref{th:ff4} that follows  concerns   the summation properties of sequences $\{a_n\}$ in $\tilde{\bf b}^{(m)}$ and  is the most important result that we use in developing the $\tilde{d}^{(m)}$ transformation. Its  proof relies
on Theorems \ref{th:ff1}, \ref{th:ff2}, and \ref{th:ff3} and  is quite involved.
We continue to use the notation of Theorem \ref{th:ff3} and $A_n=\sum^{n}_{k=1}a_k$.
\begin{theorem} [{\bf \cite{Sidi:2003:PEM}, Theorem 6.6.6}]  \label{th:ff4}
Let $\{a_n\}\in \tilde{\mathbf{b}}^{(m)}$ for which the infinite series $\sum^\infty_{n=1}a_n$  converges or diverges but has an Abel sum or Hadamard finite part.
 Then there exist a constant
$S$ and a function $g\in \tilde{\mathbf{A}}^{(0,m)}_0$ strictly
such that
\begin{equation}\label{eq:ff15}
A_{n-1}=S+n^{\sigma}a_n\, g(n),
\end{equation}
whether $\sum^{\infty}_{n=1}a_n$ converges or not.
Here $S$ is the sum of $\sum^\infty_{n=1}a_n$  when the latter converges; otherwise, $S$ is the Abel sum or the Hadamard finite part of $\sum^\infty_{n=1}a_n$.
\end{theorem}

\noindent{\bf Remark:} Before closing, we would like to mention  that we can use the $\tilde{d}^{(m)}$ transformation
 for computing the sums  of the two (trigonometric-type) series $S^{(c)}=\sum^\infty_{n=1}a^{(c)}_n$ and
$S^{(s)}=\sum^\infty_{n=1}a^{(s)}_n$ with
$$a^{(c)}_n=[\Gamma(n)]^{s/m}e^{u(n)}\cos(v(n))h(n)\quad \text{and}\quad
a^{(s)}_n=[\Gamma(n)]^{s/m}e^{u(n)}\sin(v(n))h(n),$$
where $u(n)$ and $v(n)$ are {\em real} polynomials of degree at most $m$ in $n^{1/m}$ and $h(n)\in\tilde{\bf A}_0^{(\gamma,m)}$ is not necessarily real.
Clearly, neither of the sequences $\{a^{(c)}_n\}$ or $\{a^{(s)}_n\}$ belongs to
$\tilde{\bf b}^{(m)}$. The two sequences  $\{a_n^{(\pm)}\}$, where
$$ a_n^{(\pm)}=[\Gamma(n)]^{s/m}\exp[Q^{(\pm)}(n)] h(n);\quad Q^{(\pm)}(n)=u(n)\pm\mathrm{i}v(n),$$ do belong to $\tilde{\bf b}^{(m)},$ however.
In view of this observation, we can now   apply the $\tilde{d}^{(m)}$ transformation to the two series $S^{(\pm)}=\sum^\infty_{n=1}a_n^{(\pm)}$
successfully.
Clearly,
$$ S^{(c)}=\frac{S^{(+)}+S^{(-)}}{2}\quad \text{and}\quad
S^{(s)}=\frac{S^{(+)}-S^{(-)}}{2\mathrm{i}}.$$
In case $h(n)$ is real, it is sufficient to apply the $\tilde{d}^{(m)}$ transformation to $S^{(+)}$ only since, in this case,
$$ S^{(c)}=\Re S^{(+)}\quad \text{and}\quad S^{(s)}=\Im S^{(+)}.$$

\section{The $\tilde{d}^{(m)}$ transformation} \label{se3}
\setcounter{equation}{0}
\setcounter{theorem}{0}
\subsection{Derivation of the $\tilde{d}^{(m)}$ transformation}
Consider now the cases in which  Theorem \ref{th:ff4} applies and \eqref{eq:ff15} holds.
Being in $\tilde{\mathbf{A}}^{(0,m)}_0$ strictly, the function $g(n)$ in \eqref{eq:ff15} has the asymptotic expansion
\beq \label{eq:ff25} g(n)\sim\sum^\infty_{i=0}g_in^{-i/m}\quad \text{as $n\to\infty$}, \quad g_0\neq0.\eeq
Consequently, \eqref{eq:ff15} can be expressed as in
\beq\label{eq:ff15a} A_{n-1}\sim S+n^{\sigma}a_n\sum^\infty_{i=0}g_in^{-i/m}\quad \text{as $n\to\infty$},\quad g_0\neq0.\eeq

We now go on to the development of the $\tilde{d}^{(m)}$ transformation :
First, we  truncate the infinite summation in \eqref{eq:ff15a} at the term  $i=n-1$, replace the asymptotic equality sign $\sim$ by the equality sign $=$, and replace $S$ by $\tilde{d}^{(m,j)}_n$ and the $\beta_i$ by $\bar{\beta}_i$.
Next, we choose positive integers $R_l$, $l=0,1,\ldots,$  that are ordered as in
\beq\label{eq:ff32} 1\leq R_0<R_1<R_2<\cdots, \eeq
and we set up the $(n+1)\times(n+1)$ system of linear equations
\begin{equation}\label{eq:ff26d}
A_{R_l-1}=\tilde{d}^{(m,j)}_n+R_l^{\hat{\sigma}}a_{R_l}
\sum^{n-1}_{i=0}
\frac{\bar{\beta}_i}{(R_l+\alpha)^{i/m}},\  \ j\leq l\leq j+n,
\end{equation}
where $\hat{\sigma}=\sigma$ when $\sigma$ is known or $\hat{\sigma}$ is a known upper bound for $\sigma$.
 (Needless to say, if we
know the exact value of $\sigma$, especially $\sigma=0$, we should use
it. Since $\sigma\leq m/m=1$ in all cases, we can always choose $\hat{\sigma}=1$ and be sure that the $\tilde{d}^{(m)}$ transformation will accelerate convergence successfully in all cases.)\footnote{\label{ft1} The choice $\hat{\sigma}=1$ results in  the ``universal'' formulation of the $\tilde{d}^{(m)}$ transformation that is applicable in the presence of all $\{a_n\}\in \tilde{\bf b}^{(m)}$ that we mentioned in the paragraph preceding the last in Section \ref{se1}.} Here $\alpha>-R_0$ and a good
choice in many cases  is $\alpha=0$. As can be seen from \eqref{eq:ff26d}, to compute $\tilde{d}^{(m,j)}_n$, we need the first $R_{j+n}$ terms of the infinite series, namely, $a_1,a_2,\ldots,a_{R_{j+n}}$.

Note that the unknowns in \eqref{eq:ff26} are $\tilde{d}^{(m,j)}_n$ and $\bar{\beta}_0,\bar{\beta}_1,\ldots,\bar{\beta}_{n-1}$. Of these, $\tilde{d}^{(m,j)}_n$ is the approximation to $S$ and the $\bar{\beta}_i$ are additional auxiliary unknowns. We  call this procedure  the
$\tilde{d}^{(m)}$ transformation. This transformation is actually a generalized Richardson extrapolation method in the class GREP$^{(1)}$, which is the simplest prototype of the  {\em generalized Richardson extrapolation procedure} GREP$^{(m)}$ of the author \cite{Sidi:1979:SPG}; see also Sidi \cite[Chapters 4--7]{Sidi:2003:PEM}.

The approximations $\tilde{d}^{(m,j)}_n\equiv A^{(j)}_n$ can be arranged in a two-dimensional array as in Table \ref{table1}. Note that $\tilde{d}^{(m,j)}_0=A^{(j)}_0=A_{R_j-1}$, $j=0,1,\ldots.$

\begin{table}[h!]
\renewcommand\thetable{3.1}
$$
\begin{array}{ccccc}A^{(0)}_0&&&&\\ A^{(1)}_0&A^{(0)}_1&&&\\
 A^{(2)}_0&A^{(1)}_1&A^{(0)}_2&&\\
 A^{(3)}_0&A^{(2)}_1&A^{(1)}_2&A^{(0)}_3&\\
 \vdots& \vdots& \vdots& \vdots& \ddots
\end{array}
$$
\caption{\label{table1}}
\end{table}

When $\hat{\sigma}=\hat{q}/m\geq0$, where $\hat{q}$ is a nonnegative integer, the equations in \eqref{eq:ff26d} can be replaced by
\begin{equation}\label{eq:ff26}
A_{R_l}=\tilde{d}^{(m,j)}_n+R_l^{\hat{\sigma}}a_{R_l}
\sum^{n-1}_{i=0}
\frac{\bar{\beta}_i}{(R_l+\alpha)^{i/m}},\  \ j\leq l\leq j+n,
\end{equation}
the solution for $\tilde{d}^{(m,j)}_n$ remaining the same as in \eqref{eq:ff26d}.
This amounts to adding  $a_{R_l}$ to both sides of \eqref{eq:ff26d}, and replacing $\bar{\beta}_{\hat{q}}$ by
$\bar{\beta}_{\hat{q}}+1.$
In our numerical examples, we have taken $\hat{\sigma}=m/m=1$ and used \eqref{eq:ff26} to define $\tilde{d}^{(m,j)}_n$. Note that now $A^{(j)}_0=A_{R_j}$, $j=0,1,\ldots,$ in Table \ref{table1}.
Note also that, with $\hat{\sigma}=1$ in \eqref{eq:ff26},
we do not need any further information about $Q(n)$ and the parameters $s$,  $r$, and  $\gamma$ in \eqref{eq:ff7r}; mere knowledge of the fact that $\{a_n\}$ is in $\tilde{\bf b}^{(m)}$ is sufficient for applying the  $\tilde{d}^{(m)}$ transformation successfully.

Looking at how the approximations $A^{(j)}_n$ are placed in Table \ref{table1}, we call the sequences $\{A^{(j)}_n\}_{j=0}^\infty$ (with $n$ fixed) {\em column sequences}. Similarly, we
call the sequences $\{A^{(j)}_n\}_{n=0}^\infty$ (with $j$ fixed) {\em diagonal sequences}. The known theoretical results and numerical experience suggest that  diagonal sequences have superior convergence properties and are much better than column sequences when the latter converge.
Furthermore, numerical experience suggests that diagonal sequences  converge to some antilimit when the series $\sum^\infty_{n=1}a_n$ diverges. This can be proved rigorously at least in some cases.
Normally, we look at the diagonal sequence $\{A^{(0)}_n\}_{n=0}^\infty$.

We review some of the convergence theory pertaining to the
$\tilde{d}^{(m)}$ transformation in subsection \ref{sse35}.

\subsection{Assessing the numerical stability of the $\tilde{d}^{(m)}$ transformation}
\label{sse-3.2}
An important issue that is critical at times when computing the approximations
$\tilde{d}^{(m,j)}_n$ is that of numerical stability in the presence of finite-precision arithmetic. We show how this can be tackled effectively next.

By Cramer's rule on the linear system in \eqref{eq:ff26}, $\tilde{d}^{(m,j)}_n$ can be expressed in the form
\beq  \label{eqyyy} \tilde{d}^{(m,j)}_n=\sum^n_{i=0}\gamma^{(j)}_{n,i} A_{R_{j+i}}\equiv A^{(j)}_n, \eeq
with some scalars $\gamma^{(j)}_{n,0},\gamma^{(j)}_{n,1},\ldots,\gamma^{(j)}_{n,n}$ that satisfy
$\sum^n_{i=0}\gamma^{(j)}_{n,i}=1$.
As discussed in \cite{Sidi:2003:PEM}, the numerical stability of the $A^{(j)}_n$ computed in finite-precision arithmetic can be assessed reliably as follows:
Denote the numerically computed $A_i$ and $A^{(j)}_n$ by $\bar{A}_i$ and
$\bar{A}^{(j)}_n$, respectively. Then   $\bar{A}^{(j)}_n-S$, the actual numerical error in $\bar{A}^{(j)}_n$, satisfies
\beq  |\bar{A}^{(j)}_n-S|\leq|\bar{A}^{(j)}_n-A^{(j)}_n|+|A^{(j)}_n-S|.\eeq
The term $|A^{(j)}_n-S|$ is the exact (theoretical) error in $A^{(j)}_n$ and, assuming convergence, it tends to zero as $j\to\infty$ or as $n\to\infty$. The term $|\bar{A}^{(j)}_n-A^{(j)}_n|$, however, remains  a positive quantity, meaning that the computational error  $|\bar{A}^{(j)}_n-A^{(j)}_n|$ dominates the actual  error in $\bar{A}^{(j)}_n$; that is,
\beq \label{eq:ff54}|\bar{A}^{(j)}_n-S|\approx|\bar{A}^{(j)}_n-A^{(j)}_n|\quad \text{for large $j$ or $n$}.\eeq
We now consider two different but related approaches to the estimation of $|\bar{A}^{(j)}_n-A^{(j)}_n|$, hence to the estimation of the numerical stability:

\begin{enumerate}

\item
Let us denote the absolute error in the computation of $A_i$ by $\epsilon_i$; thus $\bar{A}_i=A_i+\epsilon_i$. Then, assuming  that the computed $\gamma^{(j)}_{n,i}$ are not much different from the exact ones,\footnote{The explanation for this is twofold: (i) Numerical
computations show this. (ii) What matters is not so much the exact value
of $\Gamma^{(j)}_n$ in \eqref{eq:ff55}--\eqref{eq:ff54k} and of $\Lambda^{(j)}_n$ in
\eqref{eq:ff56}--\eqref{eq:ff59c},
 but rather
their orders of magnitude, as explained following \eqref{eq:ff59c} and as many
numerical examples show very clearly.} we have
$$\bar{A}^{(j)}_n\approx\sum^n_{i=0}\gamma^{(j)}_{n,i}\bar{A}_{R_{j+i}}=A^{(j)}_n+ \sum^n_{i=0}\gamma^{(j)}_{n,i}\epsilon_{R_{j+i}}\quad ,$$
from which, we obtain

\beq\label{eq:ff55}|\bar{A}^{(j)}_n-A^{(j)}_n|\lessapprox\Gamma^{(j)}_n \bigg(\max_{0\leq i\leq n}|\epsilon_{R_{j+i}}|\bigg), \eeq
where
\beq \label{eq:ff55a}\Gamma^{(j)}_n=\sum^n_{i=0}|\gamma^{(j)}_{n,i}|\geq 1.\eeq
Consequently, in case of convergence, \eqref{eq:ff54} becomes
\beq \label{eq:ff54a}|\bar{A}^{(j)}_n-S|\lessapprox
\Gamma^{(j)}_n \bigg(\max_{0\leq i\leq n}|\epsilon_{R_{j+i}}|\bigg), \quad\text{for large $j$ or $n$}.\eeq
If the $A_i$ are computed with machine accuracy and the roundoff unit of the floating-point arithmetic being used is ${\bf u}$, then we have $|\epsilon_i|\leq|A_i|{\bf u}$. In case the series $\sum^\infty_{n=1}a_n$ converges, we have that the $A_i$ are approximately equal to, or of the same order as,  $S$. Therefore,  \eqref{eq:ff54a} can be replaced by
\beq\label{eq:ff54k}
\frac{|\bar{A}^{(j)}_n-S|}{|S|}\lessapprox
\Gamma^{(j)}_n {\bf u}, \quad\text{for large $j$ or $n$}.\eeq
In such a case, if $\Gamma^{(j)}_n{\bf u}=O(10^{-p})$, where $p$ is a positive integer,  then the relative error in $\bar{A}^{(j)}_n$ is $O(10^{-p})$, that is, we can rely on $p$ of the significant figures of $\bar{A}^{(j)}_n$ as  being correct for $j$ or $n$ large.

Finally, by Theorem 7.2.3 in \cite[p. 161]{Sidi:2003:PEM},
\beq\label{eq:61}\Gamma^{(j)}_n=1\quad \text{if $a_na_{n+1}<0, \ \ n=1,2,\ldots.$}\eeq

\item
Let us denote the relative error in the computation of $A_i$ by $\eta_i$; thus $\bar{A}_i=A_i(1+\eta_i)$. Then, assuming again that the computed $\gamma^{(j)}_{n,i}$ are not much different from the exact ones, we have
$$\bar{A}^{(j)}_n\approx\sum^n_{i=0}\gamma^{(j)}_{n,i}\bar{A}_{R_{j+i}}=A^{(j)}_n+ \sum^n_{i=0}\gamma^{(j)}_{n,i}A_{R_{j+i}}\eta_{R_{j+i}},$$
from which, we obtain

\beq\label{eq:ff56}|\bar{A}^{(j)}_n-A^{(j)}_n|\lessapprox\Lambda^{(j)}_n \bigg(\max_{0\leq i\leq n}|\eta_{R_{j+i}}|\bigg), \eeq
where
\beq \label{eq:ff56a}\Lambda^{(j)}_n=\sum^n_{i=0}|\gamma^{(j)}_{n,i}|\, |{A}_{R_{j+i}}|.\eeq
Consequently, in case of convergence, \eqref{eq:ff54} becomes
\beq \label{eq:ff54b}|\bar{A}^{(j)}_n-S|\lessapprox
\Lambda^{(j)}_n \bigg(\max_{0\leq i\leq n}|\eta_{R_{j+i}}|\bigg), \quad\text{for large $j$ or $n$}.\eeq

The bound in \eqref{eq:ff54b} is especially useful when $\{A_n\}$ is a divergent sequence (that is, when $\sum^\infty_{n=1}a_n$ diverges) but the antilimit $S$ of $\{A_n\}$ exists and $A^{(j)}_n\to S$ as $j\to\infty$ or $n\to\infty.$\footnote{It is clear that \eqref{eq:ff54k} is useless when $\{A_n\}$ diverges.}

If the $A_i$ are computed with machine accuracy, then we have $|\eta_i|\leq{\bf u}$, where ${\bf u}$ is the roundoff unit of the floating-point arithmetic being used. In such a case, we have
\beq \label{eq:ff59}|\bar{A}^{(j)}_n-S|\lessapprox
\Lambda^{(j)}_n {\bf u}, \quad\text{for large $j$ or $n$}.\eeq
If we want to assess the relative error in $\bar{A}^{(j)}_n$, we simply divide the right-hand side of \eqref{eq:ff59} by $\bar{A}^{(j)}_n$, obtaining
\beq \label{eq:ff59c}\frac{|\bar{A}^{(j)}_n-S|}{|S|}\lessapprox
\frac{\Lambda^{(j)}_n}{|\bar{A}^{(j)}_n|} {\bf u}, \quad\text{for large $j$ or $n$,}\eeq
as an estimate of the relative error in $\bar{A}^{(j)}_n$. If
$|\Lambda^{(j)}_n/\bar{A}^{(j)}_n| {\bf u}=O(10^{-p})$ for some positive integer $p$, then we can conclude that, as an approximation to $S$,
$\bar{A}^{(j)}_n$  has approximately $p$ correct significant figures, close to convergence. Surprisingly, this seems to be the case also when the series $\sum^\infty_{n=1}a_n$   diverges weakly or strongly.
\end{enumerate}

Let us assume that the exact/theoretical diagonal sequence of approximations $\{A^{(0)}_n\}^\infty_{n=0}$
is converging to the limit or antilimit of the sequence $\{A_n\}^\infty_{n=1}$.
From our discussion above,  the following conclusion can be reached concerning the numerically computed diagonal sequence of  approximations $\{\bar{A}^{(0)}_n\}_{n=0}^\infty$: If the corresponding sequences $\{\Gamma^{(0)}_n\}_{n=0}^\infty$ and/or $\{\Lambda^{(0)}_n\}_{n=0}^\infty$ are {\em increasing} quickly, then
the accuracy of $\{\bar{A}^{(0)}_n\}_{n=0}^\infty$ is {\em decreasing} quickly, by \eqref{eq:ff54k} and/or \eqref{eq:ff59c}. Thus $\bar{A}^{(0)}_n$ may be improving (gaining more and more correct significant digits) for $n=0,1, \ldots,N,$ for some $N$, and it deteriorates for
$n=N+1,N+2,\ldots,$ in the sense that it  eventually loses all of its correct significant digits; that is, adding more terms of the series $\sum^\infty_{n=1}a_n$ in the extrapolation process does not help to improve the approximations $\bar{A}^{(0)}_n$. This is how numerical instability exhibits itself.

In subsection \ref{sse3.4}, we shall show how the  $\Gamma^{(j)}_n$ and $\Lambda^{(j)}_n$ can be computed recursively and without having to know anything other than the sequence $\{a_n\}$.

\subsection{Choice of the $R_l$}
As is obvious from \eqref{eq:ff54k} and \eqref{eq:ff59c}, the smaller $\Gamma^{(j)}_n$ and/or $\Lambda^{(j)}_n$, the better the numerical stability, hence the accuracy, of the $A^{(j)}_n$. This can be achieved by picking the integers $R_l$
in \eqref{eq:ff26d} and \eqref{eq:ff26} in one of the following two forms:
\begin{enumerate}
\item
Pick real numbers $\kappa\geq1$ and $\eta\geq1$ and set
\beq \label{aps} R_l=\lfloor \kappa l +\eta\rfloor,\quad l=0,1,\ldots.\eeq We call this choice of the $R_l$
{\em arithmetic progression sampling} and denote it by APS for short.
Clearly, $\lim_{l\to\infty}R_l/l=\kappa$, which implies that $R_l\sim \kappa l$ as $l\to\infty$, hence $\lim_{l\to\infty}R_l/R_{l-1}=1$. Note also that
\beq \kappa-1< R_l-R_{l-1}<\kappa+1\ \Rightarrow \
|(R_l-R_{l-1})-\kappa|<1,\ \forall\ l\geq1,\eeq
 whether $\kappa$ is an integer or not. Of course, the simplest APS is one in which $\kappa=1$ and $\eta=1$, that is, $R_l=l+1$, $l=0,1,\ldots.$
\item
Pick a real number $\tau>1$ and set
\beq\label{gps} R_0=1;\quad R_l=\max\{\lfloor \tau R_{l-1}\rfloor, l+1\},\quad l=1,2,\ldots.\eeq We call this choice of the $R_l$
{\em geometric progression sampling} and denote it by GPS for short.
In this case, we have (see Sidi \cite[Section 3.4]{Sidi:2012:UFE})
\beq \label{gps9} R_l=\begin{cases}l+1,&l=0,1,\ldots,L-1,\\
\lfloor \tau R_{l-1}\rfloor,& l=L,L+1,L+2,\ldots,\end{cases} \eeq
where
\beq  L=\bigg\lceil \frac{2}{\tau-1}\bigg\rceil.\eeq
In addition, $\lim_{l\to\infty}R_l/R_{l-1}=\tau$, which implies that $R_l$ increases as $\tau^l$. Indeed, GPS generates a sequence of integers $R_l$ that satisfy $b_1 \tau^l\leq R_l\leq b_2 \tau^l$ for some positive constants $b_1\leq  b_2$, hence  grow exponentially precisely like $\tau^l$.  When $\tau$ is an integer $\geq2$, then $R_l=\tau^l$ for all $l\geq0$. Of course, we do not want $R_l$ to increase very fast as this means that we need a lot of the terms of the series $\sum^\infty_{n=1}a_n$ in applying the $\tilde{d}^{(m)}$ transformation; therefore, we take $1<\tau<2$, for example.
\end{enumerate}

\noindent{\bf Remark:} Note that the sequence of the integers $R_l$ generated by APS with noninteger $\kappa$ is very closely
an arithmetic sequence, while that generated by GPS with noninteger $\tau$ is very closely a geometric sequence.

 In essentially the same form described here, APS (with integer $\kappa$ and $\eta$) and  GPS were originally suggested in Ford and Sidi \cite[Appendix B]{Ford:1987:AGR}. For a detailed discussion of the subject, see
Sidi \cite[Chapter 10]{Sidi:2003:PEM}.

\subsection{Recursive implementation via the W-algorithm} \label{sse3.4}
The W-algorithm of Sidi \cite{Sidi:1982:ASC} and its  extensions in  \cite{Sidi:1995:CAG} and  \cite[Section 7.2]{Sidi:2003:PEM} can be used to implement GREP$^{(1)}$ and study its numerical stability very efficiently. Specifically, the approximations $\tilde{d}^{(m,j)}_n$
[with $\alpha=0$ in \eqref{eq:ff26}] and the $\Gamma^{(j)}_n$ and the $\Lambda^{(j)}_n$, which are the quantities developed for assessing the numerical stability of the $\tilde{d}^{(m,j)}_n$,  can be computed very economically, {\em and without having to determine either the $\bar{\beta}_i$ in \eqref{eq:ff26} or the  $\gamma^{(j)}_{n,i}$  in \eqref{eqyyy},} as follows:
\begin{enumerate}
\item
For $j=0,1,\ldots,$ compute
\begin{align*}M^{(j)}_0&=\frac{A_{R_j}}{\omega_{R_j}},& \hspace{-1cm}
N^{(j)}_0&=\frac{1}{\omega_{R_j}}; \qquad
\omega_r\equiv r^{\hat{\sigma}}a_r,\\
H^{(j)}_0&=(-1)^j|N^{(j)}_0|,& \hspace{-1cm} K^{(j)}_0&=(-1)^j|M^{(j)}_0|.\end{align*}
\item
For $j=0,1,\ldots,$ and $n=1,2\ldots,$ compute
\begin{align*}M^{(j)}_n &= \frac{M^{(j+1)}_{n-1} - M^{(j)}_{n-1}}
{R_{j+n}^{-1/m}-R_j^{-1/m}},&\hspace{-1cm}
N^{(j)}_n &= \frac{N^{(j+1)}_{n-1} - N^{(j)}_{n-1}}
{R_{j+n}^{-1/m}-R_j^{-1/m}},\\
H^{(j)}_n&=\frac{H^{(j+1)}_{n-1} - H^{(j)}_{n-1}}
{R_{j+n}^{-1/m}-R_j^{-1/m}},&\hspace{-1cm}
K^{(j)}_n&=\frac{K^{(j+1)}_{n-1} -K^{(j)}_{n-1}}
{R_{j+n}^{-1/m}-R_j^{-1/m}}.\end{align*}
\item
For $j=0,1,\ldots,$ and $n=1,2\ldots,$ compute
$$ A^{(j)}_n=\frac{M^{(j)}_n}{N^{(j)}_n}\equiv\tilde{d}^{(m,j)}_n,
\quad \Gamma^{(j)}_n=\bigg|\frac{H^{(j)}_n}{N^{(j)}_n}\bigg|,\quad
\Lambda^{(j)}_n=\bigg|\frac{K^{(j)}_n}{N^{(j)}_n}\bigg|.$$
\end{enumerate}
Of course, the $M^{(j)}_n$, $N^{(j)}_n$, $H^{(j)}_n$, and $K^{(j)}_n$ can be arranged in separate two-dimensional tables just like the $A^{(j)}_n$ in Table \ref{table1}. For details, see \cite[Section 7.2]{Sidi:2003:PEM}.

Here we have taken $\hat{\sigma}\geq0$ and used the definition given in \eqref{eq:ff26}; hence
$M^{(j)}_0={A_{R_j}}/{\omega_{R_j}}$.  If
$\hat{\sigma}<0$, then we should use the definition given in \eqref{eq:ff26d}; therefore, $M^{(j)}_0$ should now be computed as $M^{(j)}_0={A_{R_j-1}}/{\omega_{R_j}}$.

Note that the input needed for computing $\Gamma^{(j)}_n$ and $\Lambda^{(j)}_n$ is precisely that used to compute $A^{(j)}_n$; nothing else is needed.

\subsection{Some convergence results for the $\tilde{d}^{(m)}$ transformation}
\label{sse35}

As already mentioned, the $\tilde{d}^{(m)}$ transformation is a GREP$^{(1)}$, and the convergence properties of GREP$^{(1)}$ are studied in detail in Sidi \cite{Sidi:1995:CAG}, \cite{Sidi:1999:FCS},
\cite{Sidi:2002:NCR}, and
\cite[Chapters 8, 9]{Sidi:2003:PEM}.
Powerful results on the convergence and stability of the $\tilde{d}^{(m)}$ transformation,
as it is being applied to the cases in which $\{a_n\}\in\tilde{\bf b}^{(m)}$,
can thus be found in Sidi
\cite[Chapters 8, 9]{Sidi:2003:PEM}:
\begin{itemize}
\item
For the case  $a_n=n^\gamma w(n)$, $w\in \tilde{\bf A}_0^{(0,m)}$, that is, $s=0$ and $Q(n)\equiv0$,   (mentioned in
 \cite[Example 8.2.3] {Sidi:2003:PEM}), see the  theorems in
\cite[Chapter 8]{Sidi:2003:PEM}.
\item \sloppypar
For the  cases  $a_n=e^{Q(n)}n^\gamma w(n)$ or $a_n=[\Gamma(n)]^{s/m}n^\gamma w(n)$ or
$a_n=[\Gamma(n)]^{s/m}e^{Q(n)}n^\gamma w(n)$, $w\in \tilde{\bf A}_0^{(0,m)}$, (mentioned in \cite[Example 9.2.3]{Sidi:2003:PEM}), see the  theorems in
\cite[Chapter~9]{Sidi:2003:PEM}.
\end{itemize}
Below, we state some convergence theorems that follow from those in \cite{Sidi:2003:PEM}. Here we are assuming that the functions $\mu(t)\equiv w(t^{-m})$ and $B(t)\equiv g(t^{-m})$ are both infinitely differentiable as functions of $t$ in some interval $[0,\hat{t}]$, $\hat{t}>0$. The function $g(n)$ is the one that appears in Theorem \ref{th:ff4}.

\begin{theorem}
Let $a_n=n^\gamma w(n)$ with $w\in \tilde{\bf A}_0^{(0,m)}$. Then, the following are true:
\begin{enumerate}
\item
The column sequences $\{A_n^{(j)}\}^\infty_{j=0}$  (with fixed $n$)   obtained with both  APS and GPS satisfy
$$ A_n^{(j)}-S=O(R_j\,a_{R_j}R_j^{-n/m})=O(R_j^{\gamma+1-n/m})\quad\text{as $j\to\infty$.}$$

In addition,
$ \lim_{j\to\infty}\Gamma^{(j)}_n=\infty$ {for APS} and
$\lim_{j\to\infty}\Gamma^{(j)}_n<\infty$ {for GPS}.\footnote{Thus, $\lim_{j\to\infty}A^{(j)}_n=S$ (i)\,for all $n\geq 1$ if $\sum^\infty_{k=1}a_k$
converges, that is, if $\Re\gamma<-1$, and (ii)\,for  $n>m(\Re\gamma+1)$ if $\sum^\infty_{k=1}a_k$
diverges, that is, if $\Re\gamma\geq -1$, in which case, $S$ is the antilimit.}
\item
When $\gamma$ is real, the diagonal  sequences $\{A_n^{(j)}\}^\infty_{n=0}$  (with fixed $j$) obtained with  GPS converge to $S$ whether $\sum^\infty_{k=1}a_k$ converges or not. We actually have
$$ A_n^{(j)}-S=O( e^{-\lambda n})\quad\text{as $n\to\infty$}\quad \forall\ \lambda>0.$$
This result holds also when $\gamma$ is complex, with $R_l=\tau^l$, $\tau$ being an integer.\footnote{At the present, we do not have a theorem that covers cases with  complex $\gamma$ when GPS is used with noninteger $\tau$ in \eqref{gps}.}
\end{enumerate}
\end{theorem}

\begin{theorem}
Let $a_n=e^{Q(n)}n^\gamma w(n)$, with $Q(n)=\sum^{m-1}_{i=0}\theta_i n^{1-i/m}$,
such that $\theta_0\neq0$ and  $\lim_{n\to\infty}\Re Q(n)\neq+\infty$.\footnote{This means that $a_n$ tends to zero exponentially or behaves at worst like a fixed power of $n$ as $n\to\infty$.}
 Choose $R_l$  via APS as $R_l=\kappa(l+1)$, $\kappa$ an integer. Then, the following are true:

\begin{enumerate}
\item
Provided $e^{\kappa\theta_0}\neq1$,\footnote{Note that $e^{\kappa\theta_0}=1$ only when $\theta_0$ is purely imaginary and $\kappa|\theta_0|$ is an integer multiple of $2\pi$.}
the column sequences $\{A_n^{(j)}\}^\infty_{j=0}$  (with fixed $n$) satisfy
$$ A_n^{(j)}-S=O(R_j\,a_{R_j} R_j^{-n/m} j^{-n})=
O(R_j\,a_{R_j} j^{-n/m-n})\quad\text{as $j\to\infty$}.$$
In addition,
$ \lim_{j\to\infty}\Gamma^{(j)}_n<\infty.$
 \item
 Assume $a_n$ is real and of the form $a_n=(-1)^ne^{\tilde{Q}(n)}n^\gamma w(n)$, that is,
 $Q(n)=\mrm{i}\pi n+\sum^{m-1}_{i=0}\tilde{\theta}_i n^{1-i/m}$, $\tilde{\theta}_i$ real.
 Then, whether $\sum^\infty_{k=1} a_k$ converges or not,
the diagonal  sequences $\{A_n^{(j)}\}^\infty_{n=0}$  (with fixed $j$) obtained via  APS, with $R_l=l+1$, converge to $S$. We actually have
$$ A_n^{(j)}-S=O( e^{-\lambda n})\quad\text{as $n\to\infty$}\quad \forall\ \lambda>0.$$ In addition,
$ \Gamma^{(j)}_n=1.$

 \end{enumerate}
\end{theorem}

\noindent{\bf Remarks:}
\begin{enumerate}
\item
Note that, in both theorems, $A_{R_j}-A=O(R_j^\sigma a_{R_j})$ as $j\to\infty$ by
Theorem \ref{th:ff4}; thus our results in part 1 of both theorems show clearly that convergence acceleration is taking place as $j\to\infty$ and also give  precise quantifications of the acceleration.
\item
In part 2 of both theorems,  $ A_n^{(j)}-S$ tends to zero
as $n\to\infty$ faster than {\em any}  exponential function $e^{-\lambda n}$ with  $\lambda>0$.
It is thus clear that  both theorems show that convergence acceleration is taking place as $n\to\infty$.
\end{enumerate}

\section{Illustrative examples for Theorem \ref{th:ff4}} \label{se4}
\setcounter{equation}{0}
\setcounter{theorem}{0}
We now verify Theorem \ref{th:ff4} in the form given in \eqref{eq:ff25} with a few examples, to which we will return later  in Section \ref{se5}.
The examples we choose are different kinds of telescoping series, both convergent and divergent, in which the limits or the antilimits are  identified immediately.
In these examples, we have two types of series:
\begin{align}
\text{\em Type 1:\ }&a_n=\delta_{n}-\delta_{n-1}&\Rightarrow&&
A_n&= \sum^n_{k=1}a_k=-\delta_0+\delta_{n}, \label{eq110}\\
\text{\em Type 2:\ }&a_n=(-1)^n(\delta_{n}+\delta_{n-1})&\Rightarrow&&
A_n&= \sum^n_{k=1}a_k=-\delta_0+(-1)^n\delta_{n}.\label{eq111}\end{align}

\noindent{\bf Remark.}
It seems to be quite difficult to find infinite series $\sum^\infty_{n=1}a_n$ with simple $\{a_n\}\in\tilde{\bf b}^{(m)}$ with known sums. (Our efforts to find such series in the literature have not produced any positive result.) In view of this limitation, the examples we construct here as {Type 1} and {Type 2} series are very appropriate.  As we will see shortly in Section \ref{se5}, their limits or  antilimits are simply $-\delta_0$, which are known quantities.
\medskip

\sloppypar
For simplicity, let us now  take [see  Theorem \ref{th:ff2}]
\beq
\begin{split} \label{eq112}\delta_n=(n!)^{s/m}e^{Q(n)};\quad
Q(n)=\theta_0n+\sum^{m-1}_{i=r}\theta_in^{1-i/m},\\
  s\ \ \text{integer}, \quad
\theta_0\ \ \text{real},
\quad \theta_r\neq0,\ \ r\in\{1,\ldots,m-1\}.
\end{split}
\eeq
Clearly,
\begin{itemize}
\item
when $s=0$, $\sum^\infty_{n=1}a_n$ is convergent if $\lim_{n\to\infty}\Re Q(n)=-\infty$ and divergent if $\lim_{n\to\infty}\Re Q(n)=+\infty$,
\item
when $s<0$,
$\sum^\infty_{n=1}a_n$ is always convergent, and
\item
when $s>0$, $\sum^\infty_{n=1}a_n$ is
always divergent.
\end{itemize}

We now would like to verify/confirm that, for both types of series, the sequences $\{a_n\}$ are
in $\tilde{\bf b}^{(m)}$; that is, (i)\,the relevant $a_n$ are  precisely of the form given in \eqref{eq:ff7r}--\eqref{eq:ff8r} and (ii)\,the partial sums $A_n=\sum^n_{i=1}a_i$ satisfy \eqref{eq:ff15} in Theorem \ref{th:ff4}.

\subsection{Analysis of $a_n$}
We now analyze, in a unified manner, the asymptotic behavior of $a_n$ as $n\to\infty$, recalling \eqref{eq110} and
\eqref{eq111}. By the fact that
\beq \label{eq115}\delta_{n}/\delta_{n-1}=n^{s/m}e^{f(n)};\quad f(n)=Q(n)-Q(n-1)=\Delta Q(n-1),\eeq
we have
 \begin{align} \label{eq113}
\text{Type 1:}\quad a_n&=\delta_{n-1}q(n); && q(n)=\delta_{n}/\delta_{n-1}-1=n^{s/m}e^{f(n)}-1,\\
\text{Type 2:}
\quad a_n&=(-1)^n\delta_{n-1}q(n);&& q(n)=\delta_{n}/\delta_{n-1}+1= n^{s/m}e^{f(n)}+1.\label{eq114}
\end{align}

 Now, by the binomial theorem, for $n\geq2$,
$$ \Delta(n-1)^p= n^p-(n-1)^p=n^p[1-(1-n^{-1})^p]=\sum^\infty_{j=1}(-1)^{j+1}\binom{p}{j}n^{p-j},$$
from which we have
$$ \Delta(n-1)^p=pn^{p-1}-\frac{p(p-1)}{2!}n^{p-2}+O(n^{p-3})\quad \text{as $n\to\infty$.}$$
Thus
\begin{align*}f(n)&=\theta_0+\sum^{m-1}_{i=r}\theta_i\Delta (n-1)^{1-i/m}  \notag\\ &=\theta_0+\sum^{m-1}_{i=r}\theta_i[(1-i/m)n^{-i/m}+O(n^{-1-i/m})]\quad \text{as $n\to\infty$,} \notag\\ &=\theta_0+\sum^{m-1}_{i=r}\theta_i(1-i/m)n^{-i/m}+O(n^{-1-r/m})\quad \text{as $n\to\infty$,} \notag \end{align*}
and this gives
\beq \label{eq116}f(n)=\theta_0+\sum^{m-1}_{i=r}\theta_i(1-i/m)n^{-i/m}+ u(n), \quad u(n)\in\tilde{\bf A}_0^{(-1-r/m,m)}. \eeq
Consequently,
\begin{align} \label{eq117} e^{f(n)}&=c_0\cdot\exp\bigg[\sum^{m-1}_{i=r}\theta_i(1-i/m)n^{-i/m}+ u(n)\bigg]; \quad c_0=e^{\theta_0}>0,\notag\\
&=c_0[1+ \theta_r(1-r/m)n^{-r/m}+v(n)], \quad v(n)\in\tilde{\bf A}_0^{(-(r+1)/m,m)}, \notag \\ &=c_0[1+ h(n)], \quad h(n)\in\tilde{\bf A}_0^{(-r/m,m)}\ \ \text{strictly,
because $\theta_r\neq0$}. \end{align}

We now make use of this in the analysis of $a_n$ in the two types of series:
\begin{enumerate}
\item Type 1:  By \eqref{eq113} and \eqref{eq117}, we have
$a_n=\delta_{n-1} q(n)$, where
\beq \label{eq70}q(n)= c_0n^{s/m}[1+h(n)]-1,\quad h(n)\in\tilde{\bf A}_0^{(-r/m,m)}\ \ \text{strictly}.\eeq
Then we have the following:
\beq \label{eq71} q(n)=h(n)\in\tilde{\bf A}_0^{(-r/m,m)}\ \ \text{strictly, when $s=0$ and $\theta_0=0$}. \eeq
 \beq \label{eq72}q(n)=(c_0-1)+c_0h(n)\in\tilde{\bf A}_0^{(0,m)}\ \ \text{strictly, when $s=0$ and $\theta_0\neq0$}.\eeq
 \beq \label{eq721}q(n)\in\tilde{\bf A}_0^{(s/m,m)}\ \ \text{strictly, when $s>0$}.\eeq
 \beq \label{eq722}q(n)\in\tilde{\bf A}_0^{(0,m)}\ \ \text{strictly, when $s<0$}.\eeq
\item Type 2: By \eqref{eq114}  and \eqref{eq117}, we have $a_n=(-1)^n\delta_{n-1}q(n)$, where
\beq \label{eq73}q(n)=c_0n^{s/m}[1+h(n)]+1,\quad h(n)\in\tilde{\bf A}_0^{(-r/m,m)}\ \ \text{strictly}.\eeq
Then we have the following:
\beq \label{eq74}q(n)=(c_0+1)+c_0h(n)\in\tilde{\bf A}_0^{(0,m)}\ \ \text{strictly, when $s=0$}.\eeq
\beq \label{eq75}q(n)\in\tilde{\bf A}_0^{(s/m,m)}\ \ \text{strictly, when $s>0$}.\eeq
\beq \label{eq76}q(n)\in\tilde{\bf A}_0^{(0,m)}\ \ \text{strictly, when $s<0$}.\eeq
\end{enumerate}

 Finally, let us write $a_n$ in the form
\beq \label{eq778}a_n=\delta_{n} t(n)\quad \text{for Type 1};\quad
 a_n=(-1)^n\delta_{n} t(n)\quad \text{for Type 2,}\eeq where
 \beq \label{eq779}t(n)= \frac{q(n)}{\delta_n/\delta_{n-1}}=n^{-s/m}e^{-f(n)}q(n),\eeq with the appropriate $q(n)$.
  Invoking  \eqref{eq70}--\eqref{eq76} in \eqref{eq778}-\eqref{eq779}, we conclude that
 \beq \label{eq81}\text{for Type 1:} \quad a_n=(n!)^{s/m}e^{Q(n)}n^\gamma w(n), \quad  w(n)\in\tilde{\bf A}^{(0,m)}_0\ \ \text{strictly},\eeq
 \begin{align} &\gamma=-r/m\quad \text{if $s=0$ and $\theta_0=0$};&&
 \gamma=0\quad \text{if $s=0$ and $\theta_0\neq0$}; \notag\\
 & \gamma=0\quad \text{if $s>0$;}&&
  \gamma=|s|/m\quad \text{if $s<0$.} \label{eq82}\end{align}
\beq \label{eq83}\text{for Type 2:} \quad a_n=(n!)^{s/m}e^{\widetilde{Q}(n)}n^\gamma w(n), \quad  w(n)\in\tilde{\bf A}^{(0,m)}_0\ \ \text{strictly},\eeq
\beq \label{eq84}\gamma=0\quad \text{if $s\geq0$};\quad
  \gamma=|s|/m\quad \text{if $s<0$,}\eeq
\beq \label{eq85}\widetilde{Q}(n)=\mathrm{i}\pi n+Q(n)=(\mathrm{i}\pi
+\theta_0)n+\sum^{m-1}_{i=r}\theta_in^{1-i/m}. \eeq
Here we have made use of the fact that $(-1)^n=e^{\mathrm{i}\pi n}$.

This completes the asymptotic analysis  of $a_n$ in all the different situations.
Clearly, $\{a_n\}\in \tilde{\bf b}^{(m)}$ for all the cases studied.

\subsection{Analysis of $A_n=\sum^n_{k=1}a_k$}
We now turn to the asymptotic analysis of $A_{n-1}=\sum^{n-1}_{k=1}a_k$. Let us first express   \eqref{eq110} and \eqref{eq111}, respectively,  as in

\beq A_{n-1}=-\delta_0+\delta_{n-1}\quad\text{for Type 1}\eeq
and
\beq
A_{n-1}=-\delta_0+(-1)^{n-1}\delta_{n-1}\quad\text{for Type 2}.\eeq
Invoking now (i)\,$\delta_{n-1}=a_n/q(n)$ for Type 1 and (ii)\,$\delta_{n-1}=(-1)^na_n/q(n)$ for Type 2, from \eqref{eq113} and \eqref{eq114},   respectively,  and invoking also \eqref{eq71}--\eqref{eq76},
we obtain for both Type 1 and Type 2 series

\beq \label{eq331} A_{n-1}=-\delta_0+n^\sigma a_ng(n),\quad g(n)\in\tilde{\bf A}^{(0,m)}_0\ \ \text{strictly},\eeq
with $\sigma $ assuming the following values:\\
For Type 1:
 \begin{align} &\sigma =r/m\quad \text{if $s=0$ and $\theta_0=0$};&&
 \sigma =0\quad \text{if $s=0$ and $\theta_0\neq0$}; \notag\\
 & \sigma =-s/m\quad \text{if $s>0$;}&&
  \sigma =0\quad \text{if $s<0$.} \label{eq821}\end{align}
  For Type 2:
\beq \label{eq841}\sigma=0\quad \text{if $s=0$};\quad\sigma =-s/m\quad \text{if $s>0$};\quad
  \sigma =0\quad \text{if $s<0$.}\eeq
These are clearly consistent with Theorem \ref{th:ff4} when $\delta_n=(n!)^{s/m}e^{Q(n)}$ with either $s<0$ or $s=0$ and $\lim_{n\to\infty}\Re Q(n)=-\infty$; in such cases, the series $\sum^\infty_{n=1}a_n$ converge.
Of course, Theorem \ref{th:ff4} does not directly apply to the (strongly) divergent cases for which $s>0$ or
$s=0$ and $\lim_{n\to\infty}\Re Q(n)=+\infty$, but seems to cover them too. It does so in the cases described in  \eqref{eq110}--\eqref{eq112} we have just studied.

We would like to note  that, in all the cases considered above,  $-\delta_0$ is the sum of the infinite series $\sum^\infty_{n=1}a_n$  when this series converges, that is, when $\lim_{n\to\infty}A_n$ exists; it seems to be the antilimit of $\{A_n\}$  when $\lim_{n\to\infty}A_n$ does not exist, and  numerical experiments confirm this assertion.

\section{Numerical examples I} \label{se5}
\setcounter{equation}{0}
\setcounter{theorem}{0}
We have applied
 the $\tilde{d}^{(m)}$ transformation to various infinite series $\sum^\infty_{n=1}a_n$
with $\{a_n\}\in \tilde{\bf b}^{(m)}$ for various values of $m\geq2$ and verified that it is an effective convergence accelerator.

In this section, we report  numerical results obtained from  fourteen different series with $m=2$. Specifically, we treat several   telescoping series of the types
considered in the  preceding section, namely, those with $a_n=\delta_n-\delta_{n-1}$  for Type 1 series and $a_n=(-1)^n(\delta_n+\delta_{n-1})$ for Type 2 series,   for which the limits or antilimits are known to be $S=-\delta_0$.
We also treat  examples of divergent series with unknown  antilimits.\footnote{For the divergent series considered here,  we do not even know whether antilimits exist. The approximations $A^{(0)}_n$, $n=0,1,\ldots,$  obtained by applying the $\tilde{d}^{(m)}$ transformation to these series seem definitely to converge, however. Thus, we can safely conclude that $\lim_{n\to\infty}A^{(0)}_n$ are the antilimits of these series, even though we do not know their nature.}
In our examples,  we have both (i)\,$s=0$ and $s\neq0$, (ii)\,$Q(n)\equiv0$ and  $Q(n)=\theta_0 n\pm \sqrt{n},$ (with both $\theta_0=0$ and $\theta_0\neq0$),
and (iii)\,$\gamma=0$ and $\gamma\neq0$; we observe different numerical stability issues depending on whether $s=0$ or not, $Q(n)\equiv0$ or not, and in case $Q(n)\not\equiv0$, we  observe different behavior whether $\theta_0=0$ or not.
The fact that there are several different cases, each being convergent and divergent, and each having its own convergence and stability characteristics, accounts for the large number of the numerical examples we give in this section.
Note that each example  illustrates only one of the many different cases discussed in the preceding sections.

As mentioned earlier, we can replace $\hat{\sigma}$ in \eqref{eq:ff26} by $1$, that is,
$\omega_r=ra_r$ in the W-algorithm of subsection \ref{sse3.4}, and this is what we have done here.\footnote{See footnote \ref{ft1}.} This way we do not have to worry about the exact value of $\sigma$ in \eqref{eq:ff25}. We also recall that, with $\hat{\sigma}=1$ in \eqref{eq:ff26},
we do not need any further information about $Q(n)$ and the parameters $s$,  $r$, and  $\gamma$ in \eqref{eq:ff7r}; mere knowledge of the fact that $\{a_n\}$ is in $\tilde{\bf b}^{(m)}$ is sufficient for applying the  $\tilde{d}^{(m)}$ transformation.

We have done all our computations using quadruple-precision arithmetic, for which the roundoff unit is
${\bf u}=1.93\times 10^{-34}$. This means that the highest number of significant decimal digits we can have is about $34$. In addition, if $|\Lambda^{(j)}_n/\bar{A}^{(j)}_n|=O(10^q)$ for some $q>0$, then the number of correct significant figures in $\bar{A}^{(j)}_n$ is about $p=34-q$, close to convergence.
The  tables for our numerical examples amply demonstrate the correctness of this conclusion. We advise the reader to pay attention to this fact.

In all the examples, we first try the $\tilde{d}^{(m)} $  transformation with
$R_l=l+1$, $l=0,1,\ldots,$ which is the simplest and most immediate choice for the $R_l$. As we will see, there are some slow convergence and numerical stability issues with some of these examples
when   the $R_l$  are chosen this way.
We demonstrate that these two  issues can be treated simultaneously in an effective way by using APS in some cases and GPS in some others.
In addition, it will become obvious from our numerical results that the relative error assessments shown in
\eqref{eq:ff54k} (with $\Gamma^{(j)}_n$) and in \eqref{eq:ff59c} (with $\Lambda^{(j)}_n$) are very reliable.
This clearly demonstrates the relevance and importance of the $\Gamma^{(j)}_n$ and  $\Lambda^{(j)}_n$ in assessing the accuracy of the numerical approximations to limits or antilimits. Again, the fact that
the $\Gamma^{(j)}_n$ and  $\Lambda^{(j)}_n$ can be obtained recursively via the W-algorithm, and  simultaneously with the approximations  $A^{(j)}_n$, is really surprising.

We have applied APS with
$$\text{(APS):}\quad R_l=\kappa(l+1),\quad l=0,1,\ldots; \quad\text{integer\ \ } \kappa\geq1.$$

We have applied GPS with
$$\text{(GPS):}\quad R_0=1;\quad R_l=\max\{\lfloor \tau R_{l-1}\rfloor,l+1\},\quad l=1,2,\ldots;\quad\tau\in(1,2).$$

As usual, $A_n=\sum^n_{k=1}a_k$ and $A^{(j)}_n\equiv\tilde{d}^{(m,j)}_n$,
in the tables that accompany the examples. Recall also that $R_{j+n}$ is the number of the terms of the series used for constructing $\tilde{d}^{(m,j)}_n$.

As  a rule of thumb, we can reach the following conclusions:
\begin{enumerate}
\item [C1.]
 If $a_n=h(n)\in \tilde{\bf A}_0^{(\gamma,m)}$,  use GPS.
\item [C2.]
 If $a_n=\exp[Q(n)]n^\gamma w(n)$, where $Q(n)=\sum^{m-1}_{i=1}\theta_i n^{1-i/m}$  (i.e. $\theta_0=0$), and $w(n)\in \tilde{\bf A}_0^{(0,m)}$, use GPS.
\item [C3.]
 If $a_n=\exp[Q(n)]n^\gamma w(n)$, where $Q(n)=\sum^{m-1}_{i=0}\theta_i n^{1-i/m}$  (with $\theta_0\neq0$), and $w(n)\in \tilde{\bf A}_0^{(0,m)}$, use APS.
\item [C4.]
 If   $a_n=(n!)^{s/m}\exp[Q(n)]n^\gamma w(n)$, where $s\neq0$, $Q(n)=\sum^{m-1}_{i=0}\theta_i n^{1-i/m}$, and $w(n)\in \tilde{\bf A}_0^{(0,m)}$,  use GPS.
\item [C5.]
If $a_n$ is as in any of the cases C1--C4 (with real $\theta_0$) multiplied by $(-1)^n$ for all  $n$,  use APS with $R_l=l+1$, $l=0,1,\ldots.$ Note that, in these cases, $\Gamma^{(j)}_n=1$ in accordance with \eqref{eq:61}.
 \end{enumerate}

 Of course, in all cases, we can try $R_l=l+1,$ $l=0,1,\ldots,$ first.
We use the classification  C1--C5 in our examples below.

\noindent{\bf Remark:}
Before ending, we would like to re-emphasize  the following points:
\begin{enumerate}
\item
    The only assumption we make when applying the $\tilde{d}^{(m)}$ transformation to $\sum^\infty_{n=1}a_n$ is that the sequence $\{a_n\}$ is in $\tilde{\bf b}^{(m)}$ for some $m$;
 no further information about the specific parameters [$s$, $\gamma$, $Q(n)$] in the asymptotic expansion of $a_n$ as $n\to\infty$ is needed or is being used
in the computation. We are also using the most user-friendly definition of the
$\tilde{d}^{(m)}$ transformation with $\hat{\sigma}=1$, without having to know the exact $\sigma$.
\item
 The input needed for computing $\Gamma^{(j)}_n$ and $\Lambda^{(j)}_n$ is precisely that used to compute $A^{(j)}_n$; namely, the terms $a_1,\ldots,a_{R_{j+n}}.$ Nothing else is needed. Thus all three quantities can be computed simultaneously and  efficiently by the recursive W-algorithm.
 \item
 We also recall that when
 $|\Lambda^{(j)}_n/\bar{A}^{(j)}_n| {\bf u}=O(10^{-p})$ for some positive integer $p$, we can conclude safely that, as an approximation to $S$,
$\bar{A}^{(j)}_n$  has approximately $p$ correct significant figures, close to convergence. The numbers in the  tables obtained from all of the examples below show this  to be the case both (i)\,when the series $\sum^\infty_{n=1}a_n$  converge, and also (ii)\,when they diverge, weakly or strongly.\\
To illustrate this important point, let us look at two of the (C2) examples below:
\begin{itemize}
\item
In Example  \ref{ex3}, for which the antilimit of the divergent series $\sum^\infty_{n=1}(e^{\sqrt{n}}-e^{\sqrt{n-1}})$ seems to be $S=-1$, we have the following:
In Table \ref{tab3a}, $|\Lambda^{(0)}_{28}/\bar{A}^{(0)}_{28}|{\bf u}\approx O(10^{24-34})=O(10^{-10})$,  consistent with $|\bar{A}^{(0)}_{28}-S|/|S|=O(10^{-9})$.
In Table \ref{tab3b}, $|\Lambda^{(0)}_{32}/\bar{A}^{(0)}_{32}|{\bf u}\approx O(10^{13-34})=O(10^{-21})$, consistent with $|\bar{A}^{(0)}_{32}-S|/|S|=O(10^{-19})$.
\item
In Example  \ref{ex5}, for which the antilimit of the divergent series $\sum^\infty_{n=1}e^{\sqrt{n}}$ seems to be $S=1.24628299466148185\cdots$, up to 18 decimal digits, we have the following:
In Table \ref{tab5a}, $|\Lambda^{(0)}_{28}/\bar{A}^{(0)}_{28}|{\bf u}\approx O(10^{25-34})=O(10^{-9})$,  consistent with $\bar{A}^{(0)}_{28}=1.24628299\cdots$ (the first 9 digits of $S$).
In Table \ref{tab5b}, $|\Lambda^{(0)}_{32}/\bar{A}^{(0)}_{32}|{\bf u}\approx O(10^{15-34})=O(10^{-19})$,  consistent with $\bar{A}^{(0)}_{28}=1.24628299466148185\cdots$
(the first 18 digits of $S$).
\end{itemize}
\end{enumerate}

\begin{example}\label{ex1} Let $a_n=e^{-\sqrt{n}}-e^{-\sqrt{n-1}}$, $n=1,2,\ldots.$
The series $\sum^\infty_{n=1}a_n$ is  in the C2 category and {\em converges} with limit $S=-1$.
Table \ref{tab1a} contains results obtained by choosing $R_l=l+1$, $l=0,1,\ldots.$
 In Table \ref{tab1b} we present results obtained by choosing the $R_l$ using GPS with $\tau=1.3$.
\begin{table}[h!]
\renewcommand\thetable{5.1a}
\begin{center}
$$
\mbox{$ \scriptsize
\begin{array}{||r||r|c|c|c|c||}
\hline
n&R_n& |A_{R_n}-S|/|S| &|\bar{A}^{(0)}_n-S|/|S|&\Gamma^{(0)}_n& \Lambda^{(0)}_n\\
\hline\hline
 0&      1&    3.68D-01&    3.68D-01&  1.00D+00&  6.32D-01\\
      4&      5&    1.07D-01&    7.64D-02&  2.90D+02&  2.47D+02\\
      8&      9&    4.98D-02&    4.65D-04&  3.03D+04&  2.80D+04\\
     12&     13&    2.72D-02&    1.17D-06&  4.58D+06&  4.37D+06\\
     16&     17&    1.62D-02&    1.48D-09&  7.95D+08&  7.72D+08\\
     20&     21&    1.02D-02&    1.13D-12&  1.53D+11&  1.50D+11\\
     24&     25&    6.74D-03&    5.82D-16&  3.19D+13&  3.14D+13\\
     28&     29&    4.58D-03&    4.03D-19&  7.10D+15&  7.03D+15\\
     32&     33&    3.20D-03&    2.40D-17&  1.67D+18&  1.66D+18\\
     36&    37&    2.28D-03&    9.81D-15&    4.13D+20&    4.11D+20\\
    40&    41&    1.66D-03&    1.03D-12&    1.07D+23&    1.06D+23\\
    \hline
\end{array}
$}
$$
\end{center}
\caption{\label{tab1a}  Numerical results for  Example \ref{ex1}
[$a_n=e^{-\sqrt{n}}-e^{-\sqrt{n-1}}$], where the $R_l$ are chosen using APS as $R_l=l+1$, $l=0,1,\ldots.$ Note that $S=-1$.}
\end{table}
\begin{table}[h!]
\renewcommand\thetable{5.1b}
\begin{center}
$$
\mbox{$ \scriptsize
\begin{array}{||r||r|c|c|c|c||}
\hline
n&R_n& |A_{R_n}-S|/|S| &|\bar{A}^{(0)}_n-S|/|S|&\Gamma^{(0)}_n& \Lambda^{(0)}_n\\
\hline\hline
 0&      1&    3.68D-01&    3.68D-01&  1.00D+00&  6.32D-01\\
      4&      5&    1.07D-01&    7.64D-02&  2.90D+02&  2.47D+02\\
      8&     11&    3.63D-02&    3.27D-04&  8.41D+03&  7.72D+03\\
     12&     29&    4.58D-03&    7.32D-08&  4.26D+03&  4.10D+03\\
     16&     80&    1.30D-04&    4.43D-13&  3.75D+02&  3.73D+02\\
     20&    227&    2.86D-07&    3.11D-20&  1.80D+01&  1.79D+01\\
     24&    646&    9.16D-12&    6.18D-30&  2.25D+00&  2.25D+00\\
     28&   1842&    2.29D-19&    0.00D+00&  1.10D+00&  1.10D+00\\
     32&   5258&    3.43D-32&    1.64D-33&  1.00D+00&  1.00D+00\\
    \hline
\end{array}
$}
$$
\end{center}
\vspace{-0.5cm}
\caption{\label{tab1b}  Numerical results for  Example \ref{ex1}
[$a_n=e^{-\sqrt{n}}-e^{-\sqrt{n-1}}$], where the $R_l$ are chosen using GPS with $\tau=1.3$. Note that $S=-1$.}
\end{table}
  \end{example}

\begin{example}\label{ex2} Let $a_n=(-1)^n(e^{-\sqrt{n}}+e^{-\sqrt{n-1}})$, $n=1,2,\ldots.$
The series $\sum^\infty_{n=1}a_n$ is  in the C2/C5 category and {\em converges} with limit $S=-1$.
Table \ref{tab2} contains results obtained by choosing $R_l=l+1$, $l=0,1,\ldots.$

\begin{table}[h!]
\renewcommand\thetable{5.2}
\begin{center}
$$
\mbox{$ \scriptsize
\begin{array}{||r||r|c|c|c|c||}
\hline
n&R_n& |A_{R_n}-S|/|S| &|\bar{A}^{(0)}_n-S|/|S|&\Gamma^{(0)}_n& \Lambda^{(0)}_n\\
\hline\hline
  0&      1&    3.68D-01&    3.68D-01&  1.00D+00&  1.37D+00\\
      4&      5&    1.07D-01&    4.16D-04&  1.00D+00&  1.00D+00\\
      8&      9&    4.98D-02&    3.33D-08&  1.00D+00&  1.00D+00\\
     12&     13&    2.72D-02&    6.71D-13&  1.00D+00&  1.00D+00\\
     16&     17&    1.62D-02&    5.61D-18&  1.00D+00&  1.00D+00\\
     20&     21&    1.02D-02&    2.48D-23&  1.00D+00&  1.00D+00\\
     24&     25&    6.74D-03&    6.72D-29&  1.00D+00&  1.00D+00\\
     28&     29&    4.58D-03&    4.81D-34&  1.00D+00&  1.00D+00\\
     32&     33&    3.20D-03&    3.85D-34&  1.00D+00&  1.00D+00\\
    \hline
\end{array}
$}
$$
\end{center}
\vspace{-0.5cm}
\caption{\label{tab2} Numerical results for  Example \ref{ex2}
[$a_n=(-1)^n(e^{-\sqrt{n}}+e^{-\sqrt{n-1}})$], using $R_l=l+1$, $l=0,1,\ldots.$ Note that $S=-1$.}
\end{table}
  \end{example}
\begin{example}\label{ex3} Let $a_n=e^{\sqrt{n}}-e^{\sqrt{n-1}}$, $n=1,2,\ldots.$
The series $\sum^\infty_{n=1}a_n$ is  in the C2 category and {\em diverges} with apparent antilimit $S=-1$.
Table \ref{tab3a} contains results obtained by choosing $R_l=l+1$, $l=0,1,\ldots.$
 In Table \ref{tab3b} we present results obtained by choosing the $R_l$ using GPS with $\tau=1.3$.
\begin{table}[h!]
\renewcommand\thetable{5.3a}
\begin{center}
$$
\mbox{$ \scriptsize
\begin{array}{||r||r|c|c|c|c||}
\hline
n&R_n& |A_{R_n}-S|/|S| &|\bar{A}^{(0)}_n-S|/|S|&\Gamma^{(0)}_n& \Lambda^{(0)}_n\\
\hline\hline
 0&      1&    2.72D+00&    2.72D+00&  1.00D+00&  1.72D+00\\
      4&      5&    9.36D+00&    2.60D+00&  3.13D+02&  1.65D+03\\
      8&      9&    2.01D+01&    2.10D+00&  1.03D+06&  1.17D+07\\
     12&     13&    3.68D+01&    5.08D-01&  5.04D+09&  1.01D+11\\
     16&     17&    6.18D+01&    8.73D-03&  4.72D+12&  1.51D+14\\
     20&     21&    9.78D+01&    1.64D-04&  9.85D+15&  4.75D+17\\
     24&     25&    1.48D+02&    3.52D-06&  4.11D+19&  2.86D+21\\
     28&     29&    2.18D+02&    4.14D-09&  1.23D+22&  1.20D+24\\
     32&     33&    3.12D+02&    1.67D-06&  2.04D+25&  2.73D+27\\
     36&    37&    4.38D+02&    2.14D-03&    8.15D+28&    1.46D+31\\
    40&    41&    6.04D+02&    1.52D+01&    2.76D+31&    6.56D+33\\

    \hline
\end{array}
$}
$$
\end{center}
\vspace{-0.5cm}
\caption{\label{tab3a}  Numerical results for  Example \ref{ex3}
[$a_n=e^{\sqrt{n}}-e^{\sqrt{n-1}}$], using $R_l=l+1$, $l=0,1,\ldots.$ Note that $S=-1$.}
\end{table}
\begin{table}[h!]
\renewcommand\thetable{5.3b}
\begin{center}
$$
\mbox{$ \scriptsize
\begin{array}{||r||r|c|c|c|c||}
\hline
n&R_n& |A_{R_n}-S|/|S| &|\bar{A}^{(0)}_n-S|/|S|&\Gamma^{(0)}_n& \Lambda^{(0)}_n\\
\hline\hline
 0&      1&    2.72D+00&    2.72D+00&  1.00D+00&  1.72D+00\\
      4&      5&    9.36D+00&    2.60D+00&  3.13D+02&  1.65D+03\\
      8&     11&    2.76D+01&    2.09D+00&  4.93D+05&  5.04D+06\\
     12&     29&    2.18D+02&    3.33D-01&  4.81D+07&  8.42D+08\\
     16&     80&    7.66D+03&    8.33D-04&  7.10D+07&  4.11D+09\\
     20&    227&    3.49D+06&    9.80D-08&  1.18D+07&  4.57D+09\\
     24&    646&    1.09D+11&    1.16D-11&  5.10D+06&  3.60D+10\\
     28&   1842&    4.36D+18&    3.04D-16&  1.23D+06&  5.38D+11\\
     32&   5258&    3.10D+31&    2.83D-19&  2.73D+05&  2.55D+13\\
    \hline
\end{array}
$}
$$
\end{center}
\vspace{-0.5cm}
\caption{\label{tab3b} Numerical results for  Example \ref{ex3}
[$a_n=e^{\sqrt{n}}-e^{\sqrt{n-1}}$], where the $R_l$ are chosen using GPS with $\tau=1.3$. Note that $S=-1$.}
\end{table}
  \end{example}
 \begin{example}\label{ex4} Let $a_n=(-1)^n(e^{\sqrt{n}}+e^{\sqrt{n-1}})$, $n=1,2,\ldots.$
The series $\sum^\infty_{n=1}a_n$ is  in the C2/C5 category and {\em diverges} with apparent antilimit $S=-1$.
Table \ref{tab4} contains results obtained by choosing $R_l=l+1$, $l=0,1,\ldots.$
\begin{table}[h!]
\renewcommand\thetable{5.4}
\begin{center}
$$
\mbox{$ \scriptsize
\begin{array}{||r||r|c|c|c|c||}
\hline
n&R_n& |A_{R_n}-S|/|S| &|\bar{A}^{(0)}_n-S|/|S|&\Gamma^{(0)}_n& \Lambda^{(0)}_n\\
\hline\hline
 0&      1&    2.72D+00&    2.72D+00&  1.00D+00&  3.72D+00\\
      4&      5&    9.36D+00&    1.17D-02&  1.00D+00&  5.94D+00\\
      8&      9&    2.01D+01&    4.17D-06&  1.00D+00&  1.20D+01\\
     12&     13&    3.68D+01&    2.55D-10&  1.00D+00&  2.06D+01\\
     16&     17&    6.18D+01&    5.40D-15&  1.00D+00&  3.24D+01\\
     20&     21&    9.78D+01&    5.43D-20&  1.00D+00&  4.85D+01\\
     24&     25&    1.48D+02&    3.07D-25&  1.00D+00&  6.98D+01\\
     28&     29&    2.18D+02&    1.10D-30&  1.00D+00&  9.76D+01\\
     32&     33&    3.12D+02&    3.85D-33&  1.00D+00&  1.33D+02\\
    \hline
\end{array}
$}
$$

\end{center}
\vspace{-0.5cm}
\caption{\label{tab4}  Numerical results for  Example \ref{ex4}
[$a_n=(-1)^n(e^{\sqrt{n}}+e^{\sqrt{n-1}})$], using $R_l=l+1$, $l=0,1,\ldots.$ Note that $S=-1$.}
\end{table}
  \end{example}
 \begin{example}\label{ex5} Let $a_n=e^{\sqrt{n}}$, $n=1,2,\ldots.$
The series $\sum^\infty_{n=1}a_n$ is  in the C2 category and {\em diverges}, possibly with an antilimit $S$ that is not known.
Table \ref{tab5a} contains results obtained by choosing $R_l=l+1$, $l=0,1,\ldots.$
 In Table \ref{tab5b} we present results obtained by choosing the $R_l$ using GPS with $\tau=1.3$.
\begin{table}[h!]
\renewcommand\thetable{5.5a}

\begin{center}
$$
\mbox{$ \scriptsize
\begin{array}{||r||r|c|c|c|c||}
\hline
n&R_n& A_{R_n} &\bar{A}^{(0)}_n&\Gamma^{(0)}_n& \Lambda^{(0)}_n\\
\hline\hline
 0&      1&    2.72D+00&    2.71828182845904523536028747135266231D+00&  1.00D+00&  2.72D+00\\
      4&      5&    2.92D+01&    1.25526985654591147597524366280725429D+00&  1.36D+03&  1.95D+04\\
      8&      9&    9.19D+01&    1.24677273910725248869144831373621801D+00&  3.82D+06&  1.69D+08\\
     12&     13&    2.12D+02&    1.24627554920630442266529384719068750D+00&  1.07D+10&  1.03D+12\\
     16&     17&    4.18D+02&    1.24628333168062508590422433878661035D+00&  5.25D+13&  9.49D+15\\
     20&     21&    7.52D+02&    1.24628299284466085557931505666022122D+00&  1.41D+16&  4.35D+18\\
     24&     25&    1.26D+03&    1.24628299466999557367176182613152417D+00&  1.73D+19&  8.56D+21\\
     28&     29&    2.03D+03&    1.24628299198356711403961152832133262D+00&  9.02D+22&  6.81D+25\\
     32&     33&    3.12D+03&    1.24628324241958617654595744535696746D+00&  4.54D+25&  5.06D+28\\
    \hline
\end{array}
$}
$$
\end{center}
\vspace{-0.5cm}
\caption{\label{tab5a}  Numerical results for  Example \ref{ex5}
[$a_n=e^{\sqrt{n}}$], using $R_l=l+1$, $l=0,1,\ldots.$ Note that the antilimit is not known.}
\end{table}
\begin{table}[h!]
\renewcommand\thetable{5.5b}

\begin{center}
$$
\mbox{$ \scriptsize
\begin{array}{||r||r|c|c|c|c||}
\hline
n&R_n& A_{R_n} &\bar{A}^{(0)}_n&\Gamma^{(0)}_n& \Lambda^{(0)}_n\\
\hline\hline
 0&      1&    2.72D+00&    2.71828182845904523536028747135266231D+00&  1.00D+00&  2.72D+00\\
      4&      5&    2.92D+01&    1.25526985654591147597524366280725429D+00&  1.36D+03&  1.95D+04\\
      8&     11&    1.43D+02&    1.24671505638378127944673235454440382D+00&  1.57D+06&  5.94D+07\\
     12&     29&    2.03D+03&    1.24628218609008659057856308239992091D+00&  2.73D+07&  2.12D+09\\
     16&     80&    1.26D+05&    1.24628298810418632879233471192843094D+00&  3.94D+07&  1.36D+10\\
     20&    227&    1.00D+08&    1.24628299465099003773467908406255098D+00&  1.53D+08&  5.13D+11\\
     24&    646&    5.39D+12&    1.24628299466148082969324108233096056D+00&  7.47D+06&  6.44D+11\\
     28&   1842&    3.68D+20&    1.24628299466148185365143315380862845D+00&  1.79D+06&  1.31D+13\\
     32&   5258&    4.45D+33&    1.24628299466148185195182683756923951D+00&  1.84D+06&  3.78D+15\\
    \hline
\end{array}
$}
$$
\end{center}
\vspace{-0.5cm}
\caption{\label{tab5b}  Numerical results for  Example \ref{ex5}
[$a_n=e^{\sqrt{n}}$], where the $R_l$ are chosen using GPS with $\tau=1.3$. Note that the antilimit is not known.}
\end{table}
\end{example}

\begin{example}\label{ex6} Let $a_n=(-1)^ne^{\sqrt{n}}$, $n=1,2,\ldots.$
The series $\sum^\infty_{n=1}a_n$ is  in the C2/C5 category and {\em diverges}, possibly with an antilimit $S$  that is not known.
Table \ref{tab6} contains results obtained by choosing $R_l=l+1$, $l=0,1,\ldots.$
\begin{table}[h!]
\renewcommand\thetable{5.6}

\begin{center}
$$
\mbox{$ \scriptsize
\begin{array}{||r||r|c|c|c|c||}
\hline
n&R_n& A_{R_n} &\bar{A}^{(0)}_n&\Gamma^{(0)}_n& \Lambda^{(0)}_n\\
\hline\hline
 0&      1&   -2.72D+00&   -2.71828182845904523536028747135266231D+00&  1.00D+00&  2.72D+00\\
      4&      5&   -6.22D+00&   -1.02389938149412728674621319253264616D+00&  1.00D+00&  3.42D+00\\
      8&      9&   -1.19D+01&   -1.02396073254025925488406911517938308D+00&  1.00D+00&  6.63D+00\\
     12&     13&   -2.07D+01&   -1.02396073204910424017949474942021129D+00&  1.00D+00&  1.12D+01\\
     16&     17&   -3.38D+01&   -1.02396073204906060742047437248544015D+00&  1.00D+00&  1.74D+01\\
     20&     21&   -5.26D+01&   -1.02396073204906060526543006429244242D+00&  1.00D+00&  2.58D+01\\
     24&     25&   -7.89D+01&   -1.02396073204906060526534757271439465D+00&  1.00D+00&  3.70D+01\\
     28&     29&   -1.15D+02&   -1.02396073204906060526534757003586791D+00&  1.00D+00&  5.15D+01\\
     32&     33&   -1.64D+02&   -1.02396073204906060526534757003580917D+00&  1.00D+00&  7.01D+01\\
    \hline
\end{array}
$}
$$
\end{center}
\vspace{-0.5cm}
\caption{\label{tab6}  Numerical results for  Example \ref{ex6}
[$a_n=(-1)^ne^{\sqrt{n}}$], using $R_l=l+1$, $l=0,1,\ldots.$ Note that the antilimit is not known.}
\end{table}

  \end{example}
\begin{example}\label{ex7} Let $a_n=e^{-0.2n+\sqrt{n}}-e^{-0.2(n-1)+\sqrt{n-1}}$, $n=1,2,\ldots.$
The series $\sum^\infty_{n=1}a_n$ is  in the C3 category and {\em converges} with limit $S=-1$.
Table \ref{tab7a} contains results obtained by choosing $R_l=l+1$, $l=0,1,\ldots.$
 In Table \ref{tab7b} we present results obtained by choosing the $R_l$ using APS with $\kappa=\eta=5$, thus $R_l=5(l+1)$, $l=0,1,\ldots.$
\begin{table}[h!]
\renewcommand\thetable{5.7a}

\begin{center}
$$
\mbox{$ \scriptsize
\begin{array}{||r||r|c|c|c|c||}
\hline
n&R_n& |A_{R_n}-S|/|S| &|\bar{A}^{(0)}_n-S|/|S|&\Gamma^{(0)}_n& \Lambda^{(0)}_n\\
\hline\hline
 0&     1&    2.23D+00&    2.23D+00&    1.00D+00&    1.23D+00\\
     8&     9&    3.32D+00&    3.48D+00&    1.31D+00&    3.24D+00\\
    16&    17&    2.06D+00&    3.48D+00&    9.01D+01&    1.69D+02\\
    24&    25&    1.00D+00&    3.48D+00&    1.21D+07&    1.02D+07\\
    32&    33&    4.25D-01&    3.45D+00&    4.43D+13&    8.27D+12\\
    40&    41&    1.66D-01&    2.46D-02&    8.99D+18&    3.96D+18\\
    48&    49&    6.08D-02&    3.83D-07&    1.74D+22&    1.25D+22\\
    56&    57&    2.13D-02&    5.44D-10&    3.32D+25&    2.87D+25\\
    64&    65&    7.17D-03&    9.33D-06&    6.27D+28&    5.87D+28\\
    \hline
\end{array}
$}
$$
\end{center}
\vspace{-0.5cm}
\caption{\label{tab7a}  Numerical results for  Example \ref{ex7}
[$a_n=e^{-0.2n+\sqrt{n}}-e^{-0.2(n-1)+\sqrt{n-1}}$], using $R_l=l+1$, $l=0,1,\ldots.$ Note that $S=-1$.}
\end{table}
\begin{table}[h!]
\renewcommand\thetable{5.7b}

\begin{center}
$$
\mbox{$ \scriptsize
\begin{array}{||r||r|c|c|c|c||}
\hline
n&R_n&|A_{R_n}-S|/|S| &|\bar{A}^{(0)}_n-S|/|S|&\Gamma^{(0)}_n& \Lambda^{(0)}_n\\
\hline\hline
  0&     5&    3.44D+00&    3.44D+00&    1.00D+00&    2.44D+00\\
     4&    25&    1.00D+00&    4.08D+00&    1.47D+01&    1.54D+01\\
     8&    45&    1.01D-01&    1.47D-02&    2.50D+02&    1.59D+02\\
    12&    65&    7.17D-03&    1.01D-06&    1.03D+03&    9.79D+02\\
    16&    85&    4.18D-04&    1.91D-11&    4.53D+03&    4.50D+03\\
    20&   105&    2.14D-05&    1.48D-16&    2.00D+04&    2.00D+04\\
    24&   125&    9.96D-07&    5.84D-22&    8.86D+04&    8.86D+04\\
    28&   145&    4.32D-08&    1.34D-27&    3.90D+05&    3.90D+05\\
    32&   165&    1.77D-09&    5.29D-29&    1.71D+06&    1.71D+06\\

    \hline
\end{array}
$}
$$
\end{center}
\vspace{-0.5cm}
\caption{\label{tab7b}  Numerical results for  Example \ref{ex7}
[$a_n=e^{-0.2n+\sqrt{n}}-e^{-0.2(n-1)+\sqrt{n-1}}$], where the $R_l$ are chosen using APS with $\kappa=\eta=5$; that is, $R_l=5(l+1)$, $l=0,1,\ldots.$  Note that $S=-1$.}
\end{table}
  \end{example}

\begin{example}\label{ex8} Let $a_n=(-1)^n(e^{-0.2n+\sqrt{n}}+e^{-0.2(n-1)+\sqrt{n-1}})$, $n=1,2,\ldots.$
The series $\sum^\infty_{n=1}a_n$  is  in the C3/C5 category and {\em converges} with limit $S=-1$.
Table \ref{tab8} contains results obtained by choosing $R_l=l+1$, $l=0,1,\ldots.$

\begin{table}[h!]
\renewcommand\thetable{5.8}

\begin{center}
$$
\mbox{$ \scriptsize
\begin{array}{||r||r|c|c|c|c||}
\hline
n&R_n& |A_{R_n}-S|/|S| &|\bar{A}^{(0)}_n-S|/|S|&\Gamma^{(0)}_n& \Lambda^{(0)}_n\\
\hline\hline
  0&     1&    2.23D+00&    2.23D+00&    1.00D+00&    3.23D+00\\
     4&     5&    3.44D+00&    7.30D-03&    1.00D+00&    3.12D+00\\
     8&     9&    3.32D+00&    1.43D-06&    1.00D+00&    3.45D+00\\
    12&    13&    2.73D+00&    4.81D-11&    1.00D+00&    3.24D+00\\
    16&    17&    2.06D+00&    5.60D-16&    1.00D+00&    2.80D+00\\
    20&    21&    1.47D+00&    3.09D-21&    1.00D+00&    2.30D+00\\
    24&    25&    1.00D+00&    9.60D-27&    1.00D+00&    1.82D+00\\
    28&    29&    6.60D-01&    2.06D-32&    1.00D+00&    1.40D+00\\
    32&    33&    4.25D-01&    2.12D-33&    1.00D+00&    1.12D+00\\
    \hline
\end{array}
$}
$$
\end{center}
\vspace{-0.5cm}
\caption{\label{tab8}  Numerical results for  Example \ref{ex8}
[$a_n=(-1)^n(e^{-0.2n+\sqrt{n}}+e^{-0.2(n-1)+\sqrt{n-1}})$], using $R_l=l+1$, $l=0,1,\ldots.$ Note that $S=-1$.}
\end{table}
\end{example}

\begin{example}\label{ex9} Let $a_n=e^{-0.2n+\sqrt{n}}$, $n=1,2,\ldots.$
The series $\sum^\infty_{n=1}a_n$  is  in the C3 category and {\em converges} to a limit $S$ that is not known.
Table \ref{tab9a} contains results obtained by choosing $R_l=l+1$, $l=0,1,\ldots.$
 In Table \ref{tab9b} we present results obtained by choosing the $R_l$ using APS with $\kappa=\eta=5$, thus $R_l=5(l+1)$, $l=0,1,\ldots.$
\begin{table}[h!]
\renewcommand\thetable{5.9a}
\small
\begin{center}
$$
\mbox{$ \scriptsize
\begin{array}{||r||r|c|c|c|c||}
\hline
n&R_n& A_{R_n} &\bar{A}^{(0)}_n&\Gamma^{(0)}_n& \Lambda^{(0)}_n\\
\hline\hline
  0&     1&    2.23D+00&    2.22554092849246760457953753139507683D+00&    1.00D+00&    2.23D+00\\
     8&     9&    2.85D+01&    6.91041025116709498367466955167170379D+01&    3.07D+05&    6.07D+06\\
    16&    17&    4.97D+01&    6.94975935915328024186339724785494270D+01&    6.45D+08&    2.51D+10\\
    24&    25&    6.11D+01&    6.94975997601623491994988589397430127D+01&    1.42D+12&    7.43D+13\\
    32&    33&    6.62D+01&    6.94975997602064996307484910836876701D+01&    3.03D+15&    1.83D+17\\
    40&    41&    6.83D+01&    6.94975997602064916366549830540673670D+01&    6.29D+18&    4.09D+20\\
    48&    49&    6.91D+01&    6.94975997601806743910266253659572031D+01&    1.27D+22&    8.57D+23\\
    56&    57&    6.94D+01&    6.94975997308245777099536081393349112D+01&    2.52D+25&    1.72D+27\\
    64&    65&    6.95D+01&    6.94975821438721079747682396283208104D+01&    4.89D+28&    3.37D+30\\

    \hline
\end{array}
$}
$$
\end{center}
\vspace{-0.5cm}
\caption{\label{tab9a}  Numerical results for  Example \ref{ex9}
[$a_n=e^{-0.2n+\sqrt{n}}$], using $R_l=l+1$, $l=0,1,\ldots.$ Note that the limit is not known.}
\end{table}
\begin{table}[h!]
\renewcommand\thetable{5.9b}
\small
\begin{center}
$$
\mbox{$ \scriptsize
\begin{array}{||r||r|c|c|c|c||}
\hline
n&R_n&A_{R_n} &\bar{A}^{(0)}_n&\Gamma^{(0)}_n& \Lambda^{(0)}_n\\
\hline\hline
   0&     5&    1.48D+01&    1.48469162378884426206684133887622417D+01&    1.00D+00&    1.48D+01\\
     4&    25&    6.11D+01&    6.96296394954322519778004220754137383D+01&    5.01D+01&    2.57D+03\\
     8&    45&    6.88D+01&    6.94976035518697374035719835887475438D+01&    2.02D+02&    1.34D+04\\
    12&    65&    6.95D+01&    6.94975997601454883646767731022889718D+01&    9.03D+02&    6.25D+04\\
    16&    85&    6.95D+01&    6.94975997602064989300261396945842564D+01&    4.09D+03&    2.84D+05\\
    20&   105&    6.95D+01&    6.94975997602064988540120840350043551D+01&    1.84D+04&    1.28D+06\\
    24&   125&    6.95D+01&    6.94975997602064988540031066443388019D+01&    8.24D+04&    5.72D+06\\
    28&   145&    6.95D+01&    6.94975997602064988540031067912671900D+01&    3.66D+05&    2.54D+07\\
    32&   165&    6.95D+01&    6.94975997602064988540031067892774609D+01&    1.62D+06&    1.12D+08\\
    \hline
\end{array}
$}
$$

\end{center}
\vspace{-0.5cm}
\caption{\label{tab9b}  Numerical results for  Example \ref{ex7}
[$a_n=e^{-0.2n+\sqrt{n}}$], where the $R_l$ are chosen using APS with $\kappa=\eta=5$; that is, $R_l=5(l+1)$, $l=0,1,\ldots.$  Note that the limit is not known.}
\end{table}
\end{example}
\begin{example}\label{ex10} Let $a_n=(-1)^ne^{0.2n-\sqrt{n}}$, $n=1,2,\ldots.$
The series $\sum^\infty_{n=1}a_n$   is  in the C3/C5 category and {\em diverges},   possibly with an antilimit $S$ that is not known.
Table \ref{tab10} contains results obtained by choosing $R_l=l+1$, $l=0,1,\ldots.$
\begin{table}[h!]
\renewcommand\thetable{5.10}
\begin{center}
$$
\mbox{$ \scriptsize
\begin{array}{||r||r|c|c|c|c||}
\hline
n&R_n& A_{R_n} &\bar{A}^{(0)}_n&\Gamma^{(0)}_n& \Lambda^{(0)}_n\\
\hline\hline
   0&     1&   -4.49D-01&   -4.49328964117221591430102385015562784D-01&    1.00D+00&    4.49D-01\\
     4&     5&   -3.98D-01&   -2.54755240466624695829767367307432419D-01&    1.00D+00&    2.55D-01\\
     8&     9&   -4.08D-01&   -2.54747573868605037684471045364090490D-01&    1.00D+00&    2.55D-01\\
    12&    13&   -4.43D-01&   -2.54747573734873382943870173850446560D-01&    1.00D+00&    2.55D-01\\
    16&    17&   -5.07D-01&   -2.54747573734869455944444615389169251D-01&    1.00D+00&    2.55D-01\\
    20&    21&   -6.11D-01&   -2.54747573734869455818166677569691476D-01&    1.00D+00&    2.58D-01\\
    24&    25&   -7.80D-01&   -2.54747573734869455818162487122226260D-01&    1.00D+00&    2.88D-01\\
    28&    29&   -1.05D+00&   -2.54747573734869455818162486982736643D-01&    1.00D+00&    3.62D-01\\
    32&    33&   -1.50D+00&   -2.54747573734869455818162486982731683D-01&    1.00D+00&    4.78D-01\\
    36&    37&   -2.23D+00&   -2.54747573734869455818162486982731443D-01&    1.00D+00&    6.44D-01\\
    40&    41&   -3.45D+00&   -2.54747573734869455818162486982731587D-01&    1.00D+00&    8.81D-01\\
    \hline
\end{array}
$}
$$
\end{center}
\vspace{-0.5cm}
\caption{\label{tab10}  Numerical results for  Example \ref{ex10}
[$a_n=(-1)^ne^{0.2n-\sqrt{n}}$], using $R_l=l+1$, $l=0,1,\ldots.$ Note that the antilimit is not known.}
\end{table}
\end{example}

\begin{example}\label{ex11} Let $a_n=\sqrt{n!}e^{-\sqrt{n}}-\sqrt{(n-1)!}e^{-\sqrt{n-1}}$, $n=1,2,\ldots.$
The series $\sum^\infty_{n=1}a_n$  is  in the C4 category and {\em diverges} with apparent  antilimit $S=-1$.
Table \ref{tab11a} contains results obtained by choosing $R_l=l+1$, $l=0,1,\ldots.$
 In Table \ref{tab11b} we present results obtained by choosing the $R_l$ using GPS with $\tau=1.1$.
\begin{table}[h!]
\renewcommand\thetable{5.11a}

\begin{center}
$$
\mbox{$ \scriptsize
\begin{array}{||r||r|c|c|c|c||}
\hline
n&R_n& |A_{R_n}-S|/|S| &|\bar{A}^{(0)}_n-S|/|S|&\Gamma^{(0)}_n& \Lambda^{(0)}_n\\
\hline\hline
0&     1&    3.68D-01&    3.68D-01&    1.00D+00&    6.32D-01\\
     4&     5&    1.17D+00&    3.51D-01&    1.36D+00&    7.80D-01\\
     8&     9&    3.00D+01&    3.49D-01&    2.69D+01&    3.27D+01\\
    12&    13&    2.14D+03&    3.19D-01&    4.62D+03&    7.39D+04\\
    16&    17&    3.05D+05&    7.75D-01&    5.67D+06&    1.17D+09\\
    20&    21&    7.31D+07&    7.03D-04&    4.73D+06&    1.58D+10\\
    24&    25&    2.65D+10&    4.63D-06&    4.41D+07&    2.91D+12\\
    28&    29&    1.36D+13&    3.89D-07&    7.38D+09&    1.14D+16\\
    32&    33&    9.43D+15&    9.10D-11&    4.51D+09&    1.89D+17\\
    36&    37&    8.47D+18&    4.21D-13&    6.94D+10&    8.89D+19\\
    40&    41&    9.58D+21&    1.53D-12&    4.90D+11&    2.14D+22\\
    \hline
\end{array}
$}
$$
\end{center}
\vspace{-0.5cm}
\caption{\label{tab11a}  Numerical results for  Example \ref{ex11}
[$a_n=\sqrt{n!}e^{-\sqrt{n}}-\sqrt{(n-1)!}e^{-\sqrt{n-1}}$], using $R_l=l+1$, $l=0,1,\ldots.$ Note that $S=-1$.}
\end{table}
\begin{table}[h!]
\renewcommand\thetable{5.11b}

\begin{center}
$$
\mbox{$ \scriptsize
\begin{array}{||r||r|c|c|c|c||}
\hline
n&R_n&|A_{R_n}-S|/|S| &|\bar{A}^{(0)}_n-S|/|S|&\Gamma^{(0)}_n& \Lambda^{(0)}_n\\
\hline\hline
  0&     1&    3.68D-01&    3.68D-01&    1.00D+00&    6.32D-01\\
     4&     5&    1.17D+00&    3.51D-01&    1.36D+00&    7.80D-01\\
     8&     9&    3.00D+01&    3.49D-01&    2.69D+01&    3.27D+01\\
    12&    13&    2.14D+03&    3.19D-01&    4.62D+03&    7.39D+04\\
    16&    17&    3.05D+05&    7.75D-01&    5.67D+06&    1.17D+09\\
    20&    22&    3.08D+08&    6.83D-04&    4.35D+06&    1.34D+10\\
    24&    30&    6.81D+13&    6.10D-06&    2.88D+07&    7.03D+11\\
    28&    42&    5.74D+22&    2.52D-08&    6.18D+07&    6.43D+12\\
    32&    60&    3.95D+37&    1.93D-10&    1.72D+08&    5.79D+13\\
    36&    86&    1.46D+61&    1.02D-12&    2.38D+08&    3.10D+14\\
    40&   124&    5.66D+98&    3.00D-15&    1.44D+08&    1.42D+15\\
    44&   179&    5.17D+157&    9.37D-18&    7.85D+07&    1.27D+16\\
    48&   259&    1.24D+250&    2.65D-17&    6.87D+07&    4.18D+17\\
    \hline
\end{array}
$}
$$
\end{center}
\vspace{-0.5cm}
\caption{\label{tab11b}  Numerical results for  Example \ref{ex11}
[$a_n=\sqrt{n!}e^{-\sqrt{n}}-\sqrt{(n-1)!}e^{-\sqrt{n-1}}$], where the $R_l$ are chosen using GPS with $\tau=1.1$.  Note that $S=-1$.}
\end{table}
  \end{example}
\begin{example}\label{ex12} Let $a_n=(-1)^n(\sqrt{n!}e^{-\sqrt{n}}+\sqrt{(n-1)!}e^{-\sqrt{n-1}})$, $n=1,2,\ldots.$
The series $\sum^\infty_{n=1}a_n$  is  in the C4/C5 category and {\em diverges} with apparent  antilimit $S=-1$.
Table \ref{tab12} contains results obtained by choosing $R_l=l+1$, $l=0,1,\ldots.$
\vspace{-0.5cm}
\begin{table}[h!]
\renewcommand\thetable{5.12}
\begin{center}
$$
\mbox{$ \scriptsize
\begin{array}{||r||r|c|c|c|c||}
\hline
n&R_n& |A_{R_n}-S|/|S| &|\bar{A}^{(0)}_n-S|/|S|&\Gamma^{(0)}_n& \Lambda^{(0)}_n\\
\hline\hline
 0&     1&       3.68D-01&    3.68D-01&    1.00D+00&    1.37D+00\\
     4&     5&       1.17D+00&    9.49D-04&    1.00D+00&    9.99D-01\\
     8&     9&       3.00D+01&    6.69D-07&    1.00D+00&    2.50D+00\\
    12&    13&       2.14D+03&    1.81D-10&    1.00D+00&    1.96D+01\\
    16&    17&       3.05D+05&    2.79D-14&    1.00D+00&    2.37D+02\\
    20&    21&       7.31D+07&    2.87D-18&    1.00D+00&    3.76D+03\\
    24&    25&       2.65D+10&    2.16D-22&    1.00D+00&    7.38D+04\\
    28&    29&       1.36D+13&    1.25D-26&    1.00D+00&    1.72D+06\\
    32&    33&       9.43D+15&    1.95D-28&    1.00D+00&    4.64D+07\\
    \hline
\end{array}
$}
$$
\end{center}
\vspace{-0.5cm}
\caption{\label{tab12}  Numerical results for  Example \ref{ex12}
[$a_n=(-1)^n(\sqrt{n!}e^{-\sqrt{n}}+\sqrt{(n-1)!}e^{-\sqrt{n-1}})$], using $R_l=l+1$, $l=0,1,\ldots.$ Note that $S=-1$.}
\end{table}
\end{example}
\medskip
\begin{example}\label{ex13} Let $a_n=(-1)^n\sqrt{n!}e^{-\sqrt{n}}$, $n=1,2,\ldots.$
The series $\sum^\infty_{n=1}a_n$  is  in the C4/C5 category and {\em diverges}, possibly with an antilimit $S$  that is not known.
Table \ref{tab13} contains results obtained by choosing $R_l=l+1$, $l=0,1,\ldots.$
\begin{table}[h!]
\renewcommand\thetable{5.13}

\begin{center}
$$
\mbox{$ \scriptsize
\begin{array}{||r||r|c|c|c|c||}
\hline
n&R_n& A_{R_n} &\bar{A}^{(0)}_n&\Gamma^{(0)}_n& \Lambda^{(0)}_n\\
\hline\hline
 0&     1&   -3.68D-01&   -3.67879441171442321595523770161460873D-01&    1.00D+00&    3.68D-01\\
     4&     5&   -9.65D-01&   -2.05445994186455599969353796265963566D-01&    1.00D+00&    3.18D-01\\
     8&     9&   -2.18D+01&   -2.05408671703611238491707363340804563D-01&    1.00D+00&    1.55D+00\\
    12&    13&   -1.63D+03&   -2.05408680194850300292719482233526514D-01&    1.00D+00&    1.37D+01\\
    16&    17&   -2.40D+05&   -2.05408680199979555770776037771120370D-01&    1.00D+00&    1.71D+02\\
    20&    21&   -5.88D+07&   -2.05408680199983779970682276964493250D-01&    1.00D+00&    2.80D+03\\
    24&    25&   -2.17D+10&   -2.05408680199983784675881316821330617D-01&    1.00D+00&    5.62D+04\\
    28&    29&   -1.13D+13&   -2.05408680199983784682423139886832961D-01&    1.00D+00&    1.33D+06\\
    32&    33&   -7.93D+15&   -2.05408680199983784682433513602872466D-01&    1.00D+00&    3.64D+07\\

    \hline
\end{array}
$}
$$
\end{center}
\vspace{-0.5cm}
\caption{\label{tab13}  Numerical results for  Example \ref{ex13}
[$a_n=(-1)^n\sqrt{n!}e^{-\sqrt{n}}$], using $R_l=l+1$, $l=0,1,\ldots.$ Note that the antilimit is not known.}
\end{table}

  \end{example}
\begin{example}\label{ex14} Let $a_n=n^{\sqrt{3}}/(1+\sqrt{n})$, $n=1,2,\ldots.$
The series $\sum^\infty_{n=1}a_n$  is  in the C1 category and  {\em diverges}  with an antilimit $S$ that is not known.
Table \ref{tab14a} contains results obtained by choosing $R_l=l+1$, $l=0,1,\ldots.$
 In Table \ref{tab14b} we present results obtained by choosing the $R_l$ using GPS with $\tau=1.3$. Note that $a_n=u(n)\in \tilde{\bf A}^{(\sqrt{3}-1/2,2)}_0$ and satisfies Theorem \ref{th:ff1} with $c=0$ and $b$ the antilimit in \eqref{eq:ff2}, thus $\{a_n\}\in \tilde{\bf b}^{(2)}$ by part 1 of Theorem \ref{th:ff3}.
\begin{table}[h!]
\renewcommand\thetable{5.14a}
\small
\begin{center}
$$
\mbox{$ \scriptsize
\begin{array}{||r||r|c|c|c|c||}
\hline
n&R_n& A_{R_n} &\bar{A}^{(0)}_n&\Gamma^{(0)}_n& \Lambda^{(0)}_n\\
\hline\hline
 0&     1&    5.00D-01&    5.00000000000000000000000000000000000D-01&    1.00D+00&    5.00D-01\\
     4&     5&    1.30D+01&    1.36953249947897007206016802576582650D-02&    1.05D+02&    4.99D+02\\
     8&     9&    4.83D+01&    2.14934104105471257196163245799915443D+01&    5.87D+09&    1.24D+11\\
    12&    13&    1.11D+02&   -6.36685109382204842916350986273167363D+00&    1.48D+13&    7.52D+14\\
    16&    17&    2.05D+02&   -6.33514738456051770232899035920678815D+00&    9.82D+16&    9.45D+18\\
    20&    21&    3.31D+02&   -6.33490882742039134712027408229034021D+00&    4.83D+20&    7.64D+22\\
    24&    25&    4.93D+02&   -6.33489970412111983179650131962681705D+00&    1.91D+24&    4.55D+26\\
    28&    29&    6.93D+02&   -6.33485606274876151464675840216981192D+00&    6.50D+27&    2.19D+30\\
    32&    33&    9.31D+02&   -6.23505348202801983168713408136424401D+00&    1.98D+31&    9.00D+33\\

    \hline
\end{array}
$}
$$
\end{center}
\vspace{-0.5cm}
\caption{\label{tab14a}  Numerical results for  Example \ref{ex14}
[$a_n=n^{\sqrt{3}}/(1+\sqrt{n})$], using $R_l=l+1$, $l=0,1,\ldots.$ Note that the antilimit is not known.}
\end{table}
\begin{table}[h!]
\renewcommand\thetable{5.14b}
\begin{center}
$$
\mbox{$ \scriptsize
\begin{array}{||r||r|c|c|c|c||}
\hline
n&R_n& A_{R_n} &\bar{A}^{(0)}_n&\Gamma^{(0)}_n& \Lambda^{(0)}_n\\
\hline\hline
 0&     1&    5.00D-01&    5.00000000000000000000000000000000000D-01&    1.00D+00&    5.00D-01\\
     4&     5&    1.30D+01&    1.36953249947897007206016802576582650D-02&    1.05D+02&    4.99D+02\\
     8&    11&    7.60D+01&    2.83889899046696906432360234712697271D+01&    3.39D+09&    5.96D+10\\
    12&    29&    6.93D+02&   -6.35474899562491078540752450507974049D+00&    7.70D+10&    3.14D+12\\
    16&    80&    7.04D+03&   -6.33492046280468404096402239913451429D+00&    6.28D+11&    1.26D+14\\
    20&   227&    7.53D+04&   -6.33489959483428057932844060172892160D+00&    1.50D+12&    2.25D+15\\
    24&   646&    8.00D+05&   -6.33489961177427644581715415787437553D+00&    1.99D+12&    3.00D+16\\
    28&  1842&    8.45D+06&   -6.33489961177942835930816940226182462D+00&    2.32D+12&    3.74D+17\\
    32&  5258&    8.89D+07&   -6.33489961177942862235369239450355880D+00&    2.57D+12&    4.30D+18\\
    \hline
\end{array}
$}
$$
\end{center}
\vspace{-0.5cm}
\caption{\label{tab14b}  Numerical results for  Example \ref{ex14}
[$a_n=n^{\sqrt{3}}/(1+\sqrt{n})$],
 where the $R_l$ are chosen using GPS with $\tau=1.3$. Note that the antilimit is not known.}
\end{table}
  \end{example}


\section{Application to computation of infinite products} \label{se6}
\setcounter{equation}{0}
\setcounter{theorem}{0}
The machinery of the class $\tilde{\bf b}^{(m)}$ and the $\tilde{d}^{(m)}$
transformation treated above can also be used to accelerate the convergence of some infinite products, as discussed briefly in \cite[Section 25.11]{Sidi:2003:PEM}.
Here we expand on the treatment of \cite{Sidi:2003:PEM} considerably. We deal  with  convergent infinite products\footnote{Recall that the infinite product $\prod^\infty_{n=1}f_n$ is convergent if  $\lim_{n\to\infty}\prod^n_{k=1}f_k$ exists and  is   finite and nonzero.}
of the form
\beq\label{eq:gg1} S=\prod^\infty_{n=1}(1+v_n),\quad v_n=w(n)\in\tilde{\bf A}^{(-t/m,m)}_0\ \ \text{strictly},\quad t\geq m+1 \ \ \text{integer}.\eeq
Recall that the infinite product converges if and only if $\sum^\infty_{k=1}v_k$ converges, which implies that $t/m> 1$, which in turn implies that $t\geq m+1$ since $t$ is an integer.

Let us define
\beq\label{eq:gg2} A_0=0;\quad A_n=\prod^n_{k=1}(1+v_k), \quad n=1,2,\ldots;\quad
a_n\equiv A_n-A_{n-1},\quad n=1,2,\ldots.\eeq
Then
\beq \label{eq:gg2a} A_n=\sum^n_{k=1}a_k,\quad n=1,2,\ldots,\quad \text{and}\quad S=\lim_{n\to\infty}A_n.\eeq Now,
$$ A_n=(1+v_n)A_{n-1},\quad n=2,3,\ldots.$$ Therefore,
\beq\label{eq:gg3} a_n\equiv A_n-A_{n-1}=v_nA_{n-1}\quad \Rightarrow\quad A_{n-1}=\frac{a_n}{v_n},\quad n=2,3,\ldots.\eeq
 Applying $\Delta$ to both sides of \eqref{eq:gg3}, we obtain
\beq\label{eq:gg4} a_n=\Delta A_{n-1}=\Delta(a_n/v_n)\quad\Rightarrow \quad
a_n=a_{n+1}/v_{n+1}-a_n/v_n,\eeq
 which can be written as in
\beq\label{eq:gg5} a_n=p(n)\Delta a_n;\quad p(n)=[v_{n+1}+(\Delta v_n)/v_n]^{-1}.\eeq
Now,
since $v_n\in\tilde{\bf A}^{(-t/m,m)}_0$  strictly by \eqref{eq:gg1}, we also have   $(\Delta v_n)/v_n\in\tilde{\bf A}^{(-1,m)}_0$ strictly. In addition, $-t/m<-1$. Therefore,
$1/p(n)\in\tilde{\bf A}^{(-1,m)}_0$ {strictly}, implying that $p(n)\in\tilde{\bf A}^{(1,m)}_0$ {strictly}.
This means that $\{a_n\}\in\tilde{\bf b}^{(m)}$ by Definition \ref{def:ff2}.
Consequently, the
$\tilde{d}^{(m)}$ transformation can be applied to the sequence $\{A_n\}$, hence to the series $\sum^\infty_{n=1}a_n$, successfully.

Let us now investigate the asymptotic nature of $a_n$ in more detail. We will be applying Theorem \ref{th:ff2} for this purpose.  From \eqref{eq:gg4}, we  have, in addition to \eqref{eq:gg5},
\beq\label{eq:gg6} a_{n+1}=c(n)a_n;\quad c(n)=(1+1/v_n)v_{n+1},\eeq
which, making use of the fact that  $v_{n+1}=v_n+\Delta v_n$, we can also write as
\beq\label{eq:gg7} c(n)=1+(\Delta v_n)/v_n+v_n+\Delta v_n.\eeq
Again, since   $v_n\in \tilde{\bf A}^{(-t/m,m)}_0$ strictly by \eqref{eq:gg1}, we have also
$ \Delta v_n\in \tilde{\bf A}^{(-1-t/m,m)}_0$ strictly, as a result of which, we conclude that  $(\Delta v_n)/v_n\in\tilde{\bf A}^{(-1,m)}_0$ strictly.
In addition,  since  $v_n\sim e n^{-t/m}$ as $n\to\infty$, for some constant $e\neq0$, we have
$\Delta v_n\sim e(-t/m) n^{-1-t/m}$ as $n\to\infty$, and, therefore,
$(\Delta v_n)/v_n\sim (-t/m) n^{-1}$ as $n\to\infty$. Invoking also the fact that $t\geq m+1$, we finally have that
\beq\label{eq:gg8} c(n)=1-\frac{t}{m}n^{-1}+O(n^{-1-1/m})\quad \text{as $n\to\infty$.}\eeq
Thus Theorem \ref{th:ff2} holds with $c_0=1$ and $c_1=\cdots=c_{m-1}=0$ and
$c_m=-t/m$. In addition, \eqref{eq:ff10} in Theorem \ref{th:ff2} gives
$\epsilon_1=\cdots=\epsilon_{m-1}=0$ and $\epsilon_m=-t/m$. Substituting all these into \eqref{eq:ff9}, we obtain $\mu=0$,
$\theta_0=\theta_1=\cdots=\theta_{m-1}=0$, which implies $Q(n)\equiv0$,  and $\gamma=-t/m$. As a result, \eqref{eq:ff7} gives $a_n=h(n)$, $h\in\tilde{\bf A}^{(-t/m,m)}_0$ strictly.

By  part \ref{re1} of Theorem \ref{th:ff3}, $\{A_n\}$ satisfies \eqref{eq:ff25} with $\sigma=1$. Therefore,   Theorem \ref{th:ff4} applies and we have the following result:

\begin{theorem}
Consider the convergent infinite product  $\prod^\infty_{n=1}(1+v_n)$ with $v_n=w(n)$, $w\in\tilde{\bf A}^{(-t/m,m)}_0$ strictly, $t\geq m+1$ being an integer.
Let $A_0=0$ and $A_n=\prod^n_{k=1}(1+v_k)$ and $a_n=A_n-A_{n-1}$, $n=1,2,\ldots.$ Then
\begin{equation}\label{yy1}
A_{n-1}=S+na_ng(n),\quad  \text{$g\in\tilde{\bf A}_0^{(0,m)}$\ strictly}.
\end{equation}
Therefore, we also have

\beq\label{yy2}A_n-S\sim \alpha n^{1-t/m}\quad\text{as $n\to\infty$}, \quad\text{for some $\alpha\,\neq0$}.\eeq
\end{theorem}

The asymptotic equality in \eqref{yy2} shows a very slow convergence rate for the sequence $\{A_n\}$ of the partial products in the case being considered.

Before closing, we mention that acceleration of the convergence of infinite products $\prod^\infty_{n=1}(1+v_n)$ with $\{v_n\}\in \tilde{\bf b}^{(1)}={\bf b}^{(1)}$ was first considered by   Cohen and Levin \cite{Cohen:1999:AIP}, who use a method that is in the spirit of the $d$ transformation.

\section{Numerical examples II} \label{se7}
\setcounter{equation}{0}
\setcounter{theorem}{0}
We have applied
 the $\tilde{d}^{(m)}$ transformation to  infinite products $\prod^\infty_{n=1}(1+v_n)$
with  $v_n=w(n)\in \tilde{\bf A}^{(-t/m,m)}_0$ for various values of $m\geq1$ and verified that it is an effective convergence accelerator. We discuss one example with $m=1$, for which $S$ is known, and one example with $m=2$, for which $S$ is not known.

\begin{example}\label{ex711} Let $v_n=-z^2/n^2$, $n=1,2,\ldots.$ Therefore, $m=1$.
It is known that
$ \frac{\sin\pi z}{\pi z}= \prod^\infty_{n=1}(1-\frac{z^2}{ n^2}). $
Here we show the numerical results obtained by letting $z=1/2$, for which we have
$$ S=\frac{2}{\pi}=\prod^\infty_{n=1}\bigg(1-\frac{1}{4 n^2}\bigg). $$
Table \ref{tab711a} contains results obtained by choosing $R_l=l+1$, $l=0,1,\ldots.$
 In Table \ref{tab711b} we present results obtained by choosing the $R_l$ using GPS with $\tau=1.3$.

\begin{table}[h!]
\renewcommand\thetable{7.1a}
\small
\begin{center}
$$
\mbox{$ \scriptsize
\begin{array}{||r||r|c|c|c|c||}
\hline
n&R_n&|A_{R_n}-S|/|S| &|\bar{A}^{(0)}_n-S|/|S|&\Gamma^{(0)}_n& \Lambda^{(0)}_n\\
\hline\hline
0&    1&    1.78D-01&    1.78D-01&    1.00D+00&    1.18D+00\\
  4&    5&    4.64D-02&    6.53D-03&    8.81D+01&    9.33D+01\\
  8&    9&    2.67D-02&    3.75D-06&    1.11D+04&    1.14D+04\\
 12&   13&    1.87D-02&    3.16D-10&    1.56D+06&    1.60D+06\\
 16&   17&    1.44D-02&    7.34D-15&    2.29D+08&    2.33D+08\\
 20&   21&    1.17D-02&    6.56D-20&    3.44D+10&    3.49D+10\\
 24&   25&    9.85D-03&    1.95D-21&    5.26D+12&    5.33D+12\\
 28&   29&    8.51D-03&    4.89D-19&    8.15D+14&    8.24D+14\\
 32&   33&    7.49D-03&    6.42D-17&    1.27D+17&    1.29D+17\\
    \hline
\end{array}
$}
$$

\end{center}
\vspace{-0.5cm}
\caption{\label{tab711a}  Numerical results for  Example \ref{ex711} using $R_l=l+1$, $l=0,1,\ldots.$ Note that $S=2/\pi$.}
\end{table}

\begin{table}[h!]
\renewcommand\thetable{7.1b}
\small
\begin{center}
$$
\mbox{$ \scriptsize
\begin{array}{||r||r|c|c|c|c||}
\hline
n&R_n&|A_{R_n}-S|/|S| &|\bar{A}^{(0)}_n-S|/|S|&\Gamma^{(0)}_n& \Lambda^{(0)}_n\\
\hline\hline
0&    1&    1.78D-01&    1.78D-01&    1.00D+00&    1.18D+00\\
  4&    5&    4.64D-02&    6.53D-03&    8.81D+01&    9.33D+01\\
  8&   11&    2.20D-02&    2.62D-06&    2.92D+03&    3.02D+03\\
 12&   29&    8.51D-03&    1.93D-11&    3.30D+03&    3.36D+03\\
 16&   80&    3.11D-03&    2.41D-18&    2.74D+03&    2.76D+03\\
 20&  227&    1.10D-03&    9.18D-27&    2.09D+03&    2.09D+03\\
 24&  646&    3.87D-04&    2.22D-30&    1.90D+03&    1.90D+03\\
 28& 1842&    1.36D-04&    2.38D-30&    1.82D+03&    1.82D+03\\
 32& 5258&    4.75D-05&    1.81D-29&    1.78D+03&    1.78D+03\\
    \hline
\end{array}
$}
$$
\end{center}
\vspace{-0.5cm}
\caption{\label{tab711b}  Numerical results for  Example \ref{ex711}, where the $R_l$ are chosen using GPS with $\tau=1.3$. Note that $S=2/\pi$.}
\end{table}
\end{example}


\begin{example}\label{ex712} Let $v_n=n^{-3/2}$, $n=1,2,\ldots\ .$ Therefore, $m=2$. In this case, $S$ is not known.
Table \ref{tab712a} contains results obtained by choosing $R_l=l+1$, $l=0,1,\ldots.$
 In Table \ref{tab712b} we present results obtained by choosing the $R_l$ using GPS with $\tau=1.3$.

\begin{table}[h!]
\renewcommand\thetable{7.2a}
\small
\begin{center}
$$
\mbox{$ \scriptsize
\begin{array}{||r||r|c|c|c|c||}
\hline
n&R_n& A_{R_n} &\bar{A}^{(0)}_n&\Gamma^{(0)}_n& \Lambda^{(0)}_n\\
\hline\hline
 0&    1&    2.00D+00&    2.00000000000000000000000000000000000D+00&    1.00D+00&    2.00D+00\\
  4&    5&    3.96D+00&    1.35479942920362927909756560802876967D+00&    2.68D+02&    9.27D+02\\
  8&    9&    4.82D+00&   -1.57042116160622622622411746352944801D+01&    6.43D+06&    2.80D+07\\
 12&   13&    5.35D+00&    9.53022999002732270167986334118038986D+00&    3.05D+09&    1.50D+10\\
 16&   17&    5.71D+00&    9.20517119198173662410564570330512235D+00&    3.61D+12&    1.91D+13\\
 20&   21&    5.98D+00&    9.20093271469623484626893604230623342D+00&    4.70D+15&    2.62D+16\\
 24&   25&    6.19D+00&    9.20090135744934917522373022100561298D+00&    6.32D+18&    3.68D+19\\
 28&   29&    6.37D+00&    9.20090121361950832779877985491595492D+00&    8.69D+21&    5.22D+22\\
 32&   33&    6.51D+00&    9.20090126648125570503784378346976091D+00&    1.21D+25&    7.48D+25\\
    \hline
\end{array}
$}
$$

\end{center}
\vspace{-0.5cm}
\caption{\label{tab712a}  Numerical results for  Example \ref{ex712} using $R_l=l+1$, $l=0,1,\ldots.$ Note that the limit $S$ is not known.}
\end{table}

\begin{table}[h!]
\renewcommand\thetable{7.2b}
\small
\begin{center}
$$
\mbox{$ \scriptsize
\begin{array}{||r||r|c|c|c|c||}
\hline
n&R_n& A_{R_n} &\bar{A}^{(0)}_n&\Gamma^{(0)}_n& \Lambda^{(0)}_n\\
\hline\hline
 0&    1&    2.00D+00&    2.00000000000000000000000000000000000D+00&    1.00D+00&    2.00D+00\\
  4&    5&    3.96D+00&    1.35479942920362927909756560802876967D+00&    2.68D+02&    9.27D+02\\
  8&   11&    5.11D+00&   -6.33525823410437784195106051773687585D+01&    7.84D+06&    3.33D+07\\
 12&   29&    6.37D+00&    9.25531467665752588697176832427936420D+00&    6.64D+06&    3.28D+07\\
 16&   80&    7.36D+00&    9.20093404854746290206845604435639641D+00&    1.23D+07&    7.46D+07\\
 20&  227&    8.06D+00&    9.20090121511720760645832075667746101D+00&    1.35D+07&    9.54D+07\\
 24&  646&    8.50D+00&    9.20090121315935366161482143819303662D+00&    1.37D+07&    1.08D+08\\
 28& 1842&    8.78D+00&    9.20090121315934117116570190908213652D+00&    1.42D+07&    1.19D+08\\
 32& 5258&    8.95D+00&    9.20090121315934117115672682505231045D+00&    1.45D+07&    1.26D+08\\
    \hline
\end{array}
$}
$$
\end{center}
\vspace{-0.5cm}
\caption{\label{tab712b}  Numerical results for  Example \ref{ex712}, where the $R_l$ are chosen using GPS with $\tau=1.3$. Note that the limit $S$ is not known.}
\end{table}
\end{example}


\begin{thebibliography}{10}

\bibitem{Aitken:1926:BNS}
A.C. Aitken.
\newblock On {Bernoulli's} numerical solution of algebraic equations.
\newblock {\em Proc. Roy. Soc. Edinburgh}, 46:289--305, 1926.

\bibitem{Birkhoff:1930:FTI}
G.D. Birkhoff.
\newblock Formal theory of irregular difference equations.
\newblock {\em Acta Math.}, 54:205--246, 1930.

\bibitem{Birkhoff:1932:ATS}
G.D. Birkhoff and W.J. Trjitzinsky.
\newblock Analytic theory of singular difference equations.
\newblock {\em Acta Math.}, 60:1--89, 1932.

\bibitem{Brezinski:1975:GTS}
C.~Brezinski.
\newblock G\'{e}n\'{e}ralisations de la transformation de {Shanks}, de la table
  de {Pad\'{e}}, et de l'$\epsilon$-algorithme.
\newblock {\em Calcolo}, 12:317--360, 1975.

\bibitem{Brezinski:2010:EDP}
C.~Brezinski and M.~Redivo-Zaglia.
\newblock Extensions of {Drummond's} process for convergence acceleration.
\newblock {\em Appl. Numer. Math.}, 60:1231--1241, 2010.

\bibitem{Cohen:1999:AIP}
A.M. Cohen and D.~Levin.
\newblock Accelerating infinite products.
\newblock {\em Numer. Algorithms}, 22:157--165, 1999.

\bibitem{Drummond:1976:SCT}
J.E. Drummond.
\newblock Summing a common type of slowly convergent series of positive terms.
\newblock {\em J. Austral. Math. Soc.}, Series B, 19:416--421, 1976.

\bibitem{Ford:1987:AGR}
W.F. Ford and A.~Sidi.
\newblock An algorithm for a generalization of the {Richardson} extrapolation
  process.
\newblock {\em SIAM J. Numer. Anal.}, 24:1212--1232, 1987.

\bibitem{Levin:1973:DNT}
D.~Levin.
\newblock Development of non-linear transformations for improving convergence
  of sequences.
\newblock {\em Intern. J. Computer Math.}, B3:371--388, 1973.

\bibitem{Levin:1981:TNC}
D.~Levin and A.~Sidi.
\newblock Two new classes of nonlinear transformations for accelerating the
  convergence of infinite integrals and series.
\newblock {\em Appl. Math. Comp.}, 9:175--215, 1981.
\newblock Originally appeared as a Tel Aviv University preprint in 1975.

\bibitem{Lubkin:1952:MSI}
S.~Lubkin.
\newblock A method of summing infinite series.
\newblock {\em J. Res. Nat. Bur. Standards}, 48:228--254, 1952.

\bibitem{Sablonniere:1992:ABI}
P.~Sablonni{\`e}re.
\newblock Asymptotic behaviour of iterated modified {$\Delta^2$} and
  {$\theta_2$} transforms on some slowly convergent sequences.
\newblock {\em Numer. Algorithms}, 3:401--410, 1992.

\bibitem{Sidi:1979:SPG}
A.~Sidi.
\newblock Some properties of a generalization of the {Richardson} extrapolation
  process.
\newblock {\em J. Inst. Maths. Applics.}, 24:327--346, 1979.

\bibitem{Sidi:1982:ASC}
A.~Sidi.
\newblock An algorithm for a special case of a generalization of the
  {Richardson} extrapolation process.
\newblock {\em Numer. Math.}, 38:299--307, 1982.

\bibitem{Sidi:1995:CAG}
A.~Sidi.
\newblock Convergence analysis for a generalized {Richardson} extrapolation
  process with an application to the $d^{(1)}$-transformation on convergent and
  divergent logarithmic sequences.
\newblock {\em Math. Comp.}, 64:1627--1657, 1995.

\bibitem{Sidi:1999:FCS}
A.~Sidi.
\newblock Further convergence and stability results for the generalized
  {Richardson} extrapolation process {GREP}$^{(1)}$ with an application to the
  {$D^{(1)}$}-transformation for infinite integrals.
\newblock {\em J. Comp. Appl. Math.}, 112:269--290, 1999.

\bibitem{Sidi:2002:NCR}
A.~Sidi.
\newblock New convergence results on the generalized {Richardson} extrapolation
  process {GREP}$^{(1)}$ for logarithmic sequences.
\newblock {\em Math. Comp.}, 71:1569--1596, 2002.

\bibitem{Sidi:2003:PEM}
A.~Sidi.
\newblock {\em Practical Extrapolation Methods: Theory and Applications}.
\newblock Number~10 in Cambridge Monographs on Applied and Computational
  Mathematics. Cambridge University Press, Cambridge, 2003.

\bibitem{Sidi:2012:UFE}
A.~Sidi.
\newblock A user-friendly extrapolation method for computing infinite-range
  integrals of products of oscillatory functions.
\newblock {\em IMA J. Numer. Anal.}, 32:602--631, 2012.

\bibitem{Stoer:2002:INA}
J.~Stoer and R.~Bulirsch.
\newblock {\em {Introduction to Numerical Analysis}}.
\newblock Springer-Verlag, New York, third edition, 2002.

\bibitem{VanTuyl:1976:AMA}
A.H. {Van Tuyl}.
\newblock Application of methods for acceleration of convergence to the
  calculation of singularities of transonic flows.
\newblock In {\em Pad\'{e} Approximants Method and its Applications to
  Mechanics}, number~47 in Lecture Notes in Physics, pages 209--223, Berlin,
  1976. Springer-Verlag.

\bibitem{VanTuyl:1994:ACF}
A.H. {Van Tuyl}.
\newblock Acceleration of convergence of a family of logarithmically convergent
  sequences.
\newblock {\em Math. Comp.}, 63:229--246, 1994.

\bibitem{Wimp:1974:SSW}
J.~Wimp.
\newblock The summation of series whose terms have asymptotic representations.
\newblock {\em J. Approx. Theory}, 10:185--198, 1974.

\bibitem{Wimp:1981:STA}
J.~Wimp.
\newblock {\em {Sequence Transformations and Their Applications}}.
\newblock Academic Press, New York, 1981.

\bibitem{Wynn:1956:PTN}
P.~Wynn.
\newblock On a procrustean technique for the numerical transformation of slowly
  convergent sequences and series.
\newblock {\em Proc. Cambridge Phil. Soc.}, 52:663--671, 1956.

\end{thebibliography}

\end{document}